\documentclass[12pt]{article}

\usepackage{amssymb,amsmath,amsfonts,amssymb}
\usepackage{graphics,graphicx,color,hhline}
 \usepackage{mathtools}
 \usepackage{bbm}
\usepackage[mathscr]{euscript}
\usepackage{authblk}

\def\@abssec#1{\vspace{.05in}\footnotesize \parindent .2in 
{\bf #1. }\ignorespaces} 
\graphicspath{{/EPSF/}{Figures/}}   
\providecommand{\abs}[1]{\lvert#1\rvert}
\providecommand{\norm}[1]{\lVert#1\rVert}

\newtheorem{theorem}{Theorem}[section]
\newtheorem{lemma}[theorem]{Lemma}

\def \Rm {\mathbb R}
\def \Nm {\mathbb N}
\def \Cm {\mathbb C}

\newcommand{\eps}{\varepsilon}

\newcommand{\be}{\begin{equation}}
\newcommand{\ee}{\end{equation}}
\newcommand{\bea}{\begin{eqnarray}}
\newcommand{\eea}{\end{eqnarray}}
\newcommand{\bee}{\begin{eqnarray*}}
\newcommand{\eee}{\end{eqnarray*}}

\def\fref#1{{\rm (\ref{#1})}}
\newcommand{\calQ}{\mathcal Q}

\newcommand{\calH}{\mathcal H}
\newcommand{\calL}{\mathcal L}

\newcommand{\calE}{\mathcal E}

\newcommand{\calJ}{\mathcal J}

\def\un{{\mathbbmss{1}}}

\newcommand{\TR}{\textrm{Tr} \, }

\newcommand{\cout}[1]{}

\newcommand{\rhoeq}[1]{\varrho_e[#1]}
\newcommand{\dt}{h}
\newcommand{\dx}{\Delta x}
\newcommand{\hsnorm}[1]{\norm{#1}_{\mathcal{J}_2}}
\newcommand{\jnorm}[1]{\norm{#1}_{\mathcal{J}_1}}
\newcommand{\hnorm}[1]{\norm{#1}_{\mathcal{H}}}
\newcommand{\expt}[1]{e^{-\frac{#1}{\varepsilon^2}}}

\newcommand{\tD}[1]{\widetilde{D_{#1}}}

\newcommand{\Ham}{\mathscr{H}}
\newcommand{\Vscr}{\mathscr{V}}
\newcommand{\Ascr}{\mathscr{A}}
\newcommand{\Jscr}{\mathscr{J}}
\newcommand{\Fscr}{\mathscr{F}}

\DeclareMathOperator{\Tr}{Tr}

\DeclarePairedDelimiter\bra{\langle}{\rvert}
\DeclarePairedDelimiter\ket{\lvert}{\rangle}

\begin{document}
\title{A splitting scheme for the quantum Liouville-BGK equation}
 \author{Sophia  Potoczak Bragdon \footnote{potoczak@math.colostate.edu}}
\author{Olivier Pinaud \footnote{pinaud@math.colostate.edu}}

 \affil{Department of Mathematics, Colorado State University\\ Fort Collins CO, 80523}

 \maketitle
 \begin{abstract}
   We introduce in this work an efficient numerical method for the simulation of the quantum Liouville-BGK equation, which models the diffusive transport of quantum particles. The corner stone to the model is the BGK collision operator, obtained by minimizing the quantum free energy under the constraint that the local density of particles is conserved during collisions. This leads to a large system of coupled nonlinear nonlocal PDEs whose resolution is challenging. We then define a splitting scheme that separates the transport and the collision parts, which, exploiting the local conservation of particles, leads to a fully linear collision step. The latter involves the resolution of a constrained optimization problem that is is handled with the nonlinear conjugate gradient algorithm. We prove that the time semi-discrete scheme is convergent, and as an application of our numerical scheme,  we validate the quantum drift-diffusion model that is obtained as the diffusive limit of the quantum Liouville-BGK equation.
 \end{abstract}
   
 % \abstract{We introduce in this work an efficient numerical method for the resolution of the quantum Liouville-BGK equation, which models the diffusive transport of quantum particles. The corner stone to the model is the BGK collision operator, obtained by minimizing the quantum free energy under the constraint that the local density density of particles is conserved during collisions. This leads to a large system of coupled nonlinear nonlocal PDEs}

% \begin{center}{\textbf{Abstract}}\\
% \bigskip
%   \noindent We introduce in this work an efficient numerical method for the resolution of the quantum Liouville-BGK equation, which models the diffusive transport of quantum particles. The corner stone to the model is the BGK collision operator, obtained by minimizing the quantum free energy under the constraint that the local density density of particles is conserved during collisions. This leads to a large system of coupled nonlinear nonlocal PDEs
%    \end{center}
 \section{Introduction}
 %insister maxwellian vient de Degond,  verifier $\beta$
This work is concerned with the numerical resolution of the quantum Liouville-BGK equation of the form
\begin{equation} \label{QLEintro}
    i\hbar \partial_t\varrho = [\mathscr{H},\varrho] + i\hbar \calQ(\varrho),
\end{equation}
where $\varrho$ is a density operator, i.e. a trace class self-adjoint nonnegative operator on some Hilbert space, $[\cdot,\cdot]$ denotes the commutator between two operators, $\mathscr{H}$ is a given Hamiltonian, and $\calQ$ is a BGK-type collision operator  \cite{BGK} of the form
\begin{equation*}
    \calQ(\varrho) = \frac{1}{\tau}(\rhoeq{\varrho}-\varrho).
\end{equation*}
Above, $\tau$ is a given relaxation time, and $\rhoeq{\varrho}$ is a quantum statistical equilibrium that will be discussed further. This problem is motivated by a series of papers by Degond and Ringhofer on the derivation of quantum hydrodynamical models from first principles. In \cite{DR}, their main idea is to transpose to the quantum setting the entropy closure strategy that Levermore used for kinetic equations \cite{levermore}. As in the kinetic case, an infinite cascade of equations for the local moments of $\varrho$ can be derived from \fref{QLEintro}, and this cascade cannot be closed since  moments of a given order depend on moments of higher order. The local moments of $\varrho$ can be defined in terms of the Wigner transform  $W(x,p)$ of $\varrho$, see e.g. \cite{LP}, and by computing moments with respect to $p$,  yielding then functions of the spatial variable $x$ such as the local density of particles, the local current, and the local energy. By analogy with the classical case, Degond and Ringhofer then introduce a quantum statistical equilibrium $\varrho_e[\varrho]$ that is used to close the moments hierarchy. Depending on the number of moments accounted for in the closure procedure, several quantum macroscopic models can be obtained: Quantum Euler, Quantum Energy Transport, Quantum Navier-Stokes, or Quantum Drift-Diffusion in the diffusive regime, we refer to \cite{DGMRlivre,QHD-CMS,isotherme,QHD-review,QET,jungel-matthes-milisic,jungel-matthes,jungelbook} for more details about these models and other references on quantum hydrodynamics. The quantum Liouville-BGK equation is the ``mother'' of all of these quantum hydrodynamical models, and is therefore an important equation to study.

As in e.g. \cite{DGMRlivre,QET}, we consider in this work the case where $\varrho_e[\varrho]$ is obtained by minimizing the quantum free energy $\Fscr$, which is defined by, for appropriate density operators $\sigma$, 
\begin{equation}
\label{FE}
    \Fscr(\sigma) = k_{B}T_{0}\TR(\sigma \log \sigma - \sigma) + \TR(\Ham \sigma),
\end{equation}
under the constraint that the local density of particles of $\sigma$ is the same as that of $\varrho$, where $\varrho\equiv \varrho(t)$ is the solution to \fref{QLEintro}. In other words, if $W_\sigma$ and $W_\varrho$ are the Wigner transforms of $\sigma$ and $\varrho$, then this constraint is expressed mathematically as
$$n_\sigma:=\int W_\sigma(x,p) dp=\int W_\varrho(x,p) dp=n_\varrho,$$
(we will use more a convenient form for the definition of $n_\sigma$ later). In \fref{FE}, $\Tr(\cdot)$ denotes operator trace, $k_{B}$ is the Boltzmann constant, and $T_{0}$ is the temperature. This model gives rise, in the diffusion limit valid at time scales much larger than $\tau$, to the Quantum Drift-Diffusion model, see e.g. \cite{isotherme}. The latter is a generalization of the classical drift-diffusion model that accounts for quantum effects in a non-perturbative manner.

At the mathematical level, \fref{QLEintro} is studied in \cite{MP-JDE} in a one-dimensional spatial domain, and the minimization of $\Fscr$ under various configurations is addressed in \cite{MP-JSP,MP-KRM,DP-JMPA,DP-JFA,DP-CVPDE,DP-CMP,DP-MB}. Note also that the equilibrium $\varrho_e[\varrho]$ is central in the work of Nachtergale and Yau in their derivation of the Euler equations of fluid dynamics from many-body quantum mechanics, see \cite{nachter}.

Our main motivation in this work is to develop an efficient numerical method for the resolution of the quantum Liouville-BGK equation \fref{QLEintro}. As the minimizer of $\Fscr$ under the density constraint, the equilibrium operator $\varrho_e[\varrho]$ depends nonlinearly and nonlocally on $\varrho$, and \fref{QLEintro} can then be seen as an infinite system of coupled nonlinear nonlocal PDEs. The main difficulty in the calculation is naturally to properly handle $\varrho_e[\varrho]$. We propose here a simple and effective way to proceed by using a splitting scheme, and treat the transport term $[\Ham,\varrho]$ and the collision term $i \hbar \calQ(\varrho)$ separately. The key point is, by construction, that the local density $n_\varrho$ is a collision invariant, and as a consequence the solution to the collision step
$$  \partial_t \varrho=\calQ(\varrho), \qquad \varrho(t_0)=\sigma,
$$
satisfies $n_{\varrho(t)}=n_{\sigma}$ for all $t \geq t_0$. This yields that $\calQ(\varrho)$ takes the form $\calQ(\varrho)=\varrho_e[\sigma]-\varrho$ and is linear in $\varrho$. There is still a constrained optimization problem to solve at each time step to obtain $\varrho_e[\sigma]$, but the originally nonlinear problem is now linear. The treatment of $[\Ham,\varrho]$ is standard and poses no particular difficulty.

While the method generalizes immediately to two and three dimensional spatial settings, we will for simplicity implement and study this splitting scheme in a one-dimensional framework. One-dimensional models are revelant for instance in the study of quantum heterostructures formed by stacking layers of different materials along one direction, here $x$: electrons in the conduction band see sharp changes in the potential along $x$, while variations are small in the transverse plane; the transport properties in the bulk of the material are then calculated by imposing periodic boundary conditions in the transverse plane. We will explain informally how to derive a 1D model from a 3D one in this context in the Appendix.

Our main contributions in this work are the following: (i) implementation and analysis of a splitting scheme for \fref{QLEintro}; we will prove that the splitting solution converges to the original solution, and a by-product of the proof is the uniqueness of solutions to \fref{QLEintro} while only existence was obtained in \cite{MP-JDE}; (ii) as an application of the numerical method, we validate the Quantum Drift-Diffusion model (QDD) defined further; we compare the solutions to \fref{QLEintro} for various collision strengths with those of QDD and show an excellent agreement in the regime of validity of QDD.% (iii) we compare the solutions to \fref{QLEintro} with a non collisional model such as the Schr\"odinger-Poisson system.

The paper is structured as follows: we define in detail in Section \ref{models} the quantum Liouville-BGK equation and its diffusive limit, the QDD model. We present in Section \ref{numQLE} our numerical method for the resolution of the quantum Liouville equation: we introduce the temporal and spatial discretizations, and show that the unique time-discrete solution given by a Strang splitting scheme converges to the unique solution to the Liouville equation. The resolution of the QDD model is addressed in Section \ref{numQDD}. The numerical simulations and some algorithmic details are offered in Section \ref{simu}. Finally, an Appendix collects various technical results needed throughout the article.

\paragraph{Acknowledgment.} This work is supported by NSF CAREER Grant DMS-1452349 and NSF grant DMS-2006416. 
\section{Models} \label{models}

We introduce in this section the Quantum Liouville-BGK equation and the Quantum Drift-Diffusion model.

\subsection{The Quantum Liouville-BGK Equation (QLE)} \label{subQLE}
 We first write a density operator $\varrho$ in terms of its spectral elements,
\[\varrho = \sum_{p\in \mathbb{N}} \rho_{p}\ket{\psi_{p}}\bra{\psi_{p}},\]
where we used the Dirac bra-ket notation, and where $\{\rho_p,\psi_p\}$ are the $p-$th eigenvalue and eigenfunction pair for $\varrho$, eigenvalues counted with multiplicity. In our problem of interest, the density operators are typically full-rank, that is all eigenvalues are strictly positive, and form then a sequence $\{\rho_p\}_{p \in \Nm}$ decreasing to zero. This is a consequence of the fact, proved in \cite{MP-JSP}, that the equilibrium $\varrho_e[\varrho]$ is full-rank. With this notation, the local density $n_\varrho$ associated to $\varrho$ is defined by
\[n_{\varrho} := \sum_{p\in \mathbb{N}} \rho_{p}\abs{\psi_{p}}^2.\]
The local density can also be equivalently defined by duality in terms of the trace operator  $\Tr(\cdot)$, i.e., with $[0,L]$ our spatial domain, 

 \[(n_{\varrho},\psi) := \int_{[0,L]} n_{\varrho}\psi dx = \Tr(\varrho\psi),\] for all smooth function $\psi$ (we identify $\psi$ with the corresponding multiplication operator).

In the context of particle transport in nanostructures, the Hamiltonian $\Ham$ in \fref{QLEintro} is given by 
\[  \Ham = \Ham_0 -e\Vscr^{\text{ext}}-e\Vscr,  \quad \text{ with } \quad\Ham_0 = -\frac{\hbar^2}{2m^*}\Delta,\]
where $\Delta=d^2/dx^2$, $m^*$ is the effective mass of the electron (assumed for simplicity to be constant in the domain; considering a varying $m^*$ would only require minor modifications), and $-e$ is the electron charge. In $\Ham$, $\Vscr^{\text{ext}}$ is a bounded externally applied potential, and $\Vscr$ is electrostatic potential solution the Poisson equation 
%\begin{equation}
%\label{poisson}
$$
    \epsilon_0 \Delta \Vscr = n_\varrho, \quad \Vscr(0)=\Vscr(L)=0.
$$
%\end{equation}
Above, $\epsilon_0$ is the permittivity of the material (assumed once more to be constant for simplicity), and the maximum principle shows that $\Vscr$ is negative. The Hamiltonians  $\Ham$ and $\Ham_0$ are equipped with Neumann boundary conditions and are defined on the following domain
\be \label{domain} D(\Ham)=D(\Ham_0) = \left\{\varphi \in H^2(0,L): \frac{d}{dx} \varphi(0) = \frac{d}{dx}\varphi(L) = 0 \right\},\ee
where $H^2(0,L)$ is the usual Sobolev space. With such boundary conditions, the total number of particles in the system is fixed, and there is no particle current at the boundary. We will then model the inflow of particles by using superpositions of wave packets located away from the boundary as initial conditions. A better way to include particle flow into the domain is to use transparent boundary conditions as e.g. in \cite{ARNOLD-VLSI,PINAUD-JAP}, but this is quite technical and beyond the scope of this work. Neumann boundary conditions are chosen over homogeneous Dirichlet boundary conditions since they ensure that the density $n_\varrho$ is strictly positive over the domain. Spatial points where $n_\varrho$ vanishes (i.e. where there is no particle) are problematic when solving the minimization problem, and are then avoided with Neumann conditions, see e.g. \cite{MP-JDE} for a discussion of this matter.

Regarding the calculation of the equilibrium and the minimization of the free energy, it is shown formally in \cite{DR,isotherme} (and rigorously in \cite{MP-JSP,DP-JFA}), that $\varrho_e[\varrho]$ takes on the form of a so-called ``quantum Maxwellian'',
\begin{equation} \label{quantM}
    \rhoeq{\varrho} = e^{-(\Ham_0+\Ascr[\varrho])/k_{B}T_{0}},
\end{equation}
where $\Ascr[\varrho](t,x)$ is the chemical potential obtained as the Lagrange multiplier associated with the local density constraint $n_\sigma=n_\varrho$. It is moreover shown in \cite{DR} that the constrained optimization problem can be reformulated as the unconstrained minimization of the following convex functional of $\Ascr$:
\be \label{defJ}
\Jscr(\Ascr) = k_{B}T_{0}\, \TR\left(e^{-(\Ham_0+\Ascr)/k_{B}T_{0}}\right) + \int_{[0,L]}n_\varrho\, \Ascr\, dx.\ee
As for the density operator, we can represent the quantum Maxwellian in terms of the spectral elements $\{\lambda_{p}[\Ascr_\varrho],\phi_{p}[\Ascr_\varrho]\}_{p \in \Nm}$ of the Hamiltonian $\Ham_{\Ascr} = \Ham_{0} + \Ascr[\varrho]$, $\Ascr_\varrho\equiv \Ascr[\varrho]$ so we have
\[\rhoeq{\varrho} = \sum_{p\in \mathbb{N}} e^{-\lambda_{p}[\Ascr_\varrho]/k_{B}T_{0}}\ket{\phi_{p}[\Ascr_\varrho]}\,\bra{\phi_{p}[\Ascr_\varrho]}.\]

Following the scalings used in \cite{QDDscale}, we nondimentionalize QLE in a manner that incorporates the relevant physical constants. The characteristic length is determined by the size of the device, $\bar{x} = L$; the relaxation time is $\tau=\frac{m^* \mu}{e}$, where  $\mu$ is the (supposed constant) mobility of the electrons in the material; the reference time is given by $\bar{t} = \tfrac{L^2e}{\mu k_{B}T_{0}}$; voltages are scaled with respect to the thermal potential $\bar{V} = \tfrac{k_{B}T_{0}}{e}$, and densities with respect to the uniform density $\bar n=L^{-1}$. Using these reference values, we can now define the following dimensionless quantities:
\begin{equation}
\label{scaling}
    x' = \frac{x}{\bar{x}}, \quad n' = \frac{n}{\bar{n}}, \quad t' = \frac{t}{\bar{t}}, \quad \Vscr' = \frac{\Vscr}{\bar{V}}, \quad \Ascr' = \frac{\Ascr}{e \bar{V}},
\end{equation}
to obtain the scaled QLE coupled with the Poisson equation (omitting the primes):
\begin{equation*}
\left\{
\begin{split}
&i\varepsilon \partial_t\varrho = \frac{1}{\sqrt{2}\beta}[\Ham,\varrho] + \frac{i}{\varepsilon}(\rhoeq{\varrho} -\varrho), \quad x\in (0,1) \\
&\alpha^2\Delta \Vscr = n_\varrho, \quad \Vscr(0)=\Vscr(1)=0.
\end{split}\right.
\end{equation*}
Above, the Hamiltonian is given by 
\begin{equation*}
    \Ham = -\beta^2\Delta -\Vscr- \Vscr^{\text{ext}} =: \Ham_0 -\Vscr- \Vscr^{\text{ext}}.
\end{equation*}
The equilibrium operator $\varrho_{e}$ is 
\begin{equation*}
    \rhoeq{\varrho} = e^{-(\Ham_0+\Ascr[\varrho])},
\end{equation*}
and the dimensionless constants are
\begin{equation*}
%\label{dimlessconsts}
    \alpha = \sqrt{\frac{\epsilon_0 k_{B}T_{0}}{e^2L^2\bar{n}}}=\frac{\lambda_{d}}{L}, \quad \beta = \sqrt{\frac{\hbar^2}{2m^{*}L^2k_{B}T_{0}}}= \frac{\lambda_{dB}}{L}, \quad \varepsilon = \sqrt{\frac{k_{B}T_{0}\tau^2}{m^{*}L^2}} = \frac{\lambda_{mfp}}{L},
\end{equation*}
where $\lambda_{d}$ is the Debye length, $\lambda_{dB}$ is the de Broglie length, and $\lambda_{mfp}$ is the mean free path. We will consider moderate values $\eps=0.1$ to small values of $\eps=0.0025$ to validate the QDD model. The parameter $\beta$ controls the oscillations in the solution. Interesting (and more computationally involved) regimes correspond to small $\beta$, where particles travel large distances in the device and have wavelengths comparable with variations in the potentials. Note that small values of $\beta$ allow for a significant number of modes in the quantum Maxwellian, which justifies the use of mixed states. The parameter $\alpha$ has a relatively weak influence on the solutions. %\tr{typical values for $\eps$, $\alpha$ and $\beta$ in a given material at a given temperature}

We now turn to the Quantum Drift-Diffusion model.

\subsection{The Quantum Drift-Diffusion Model (QDD)}
QDD is obtained as the diffusive limit of QLE, i.e. in the limit as $\eps \to 0$, see \cite{QDDscale} for a derivation. The dimensional quantities in QDD  are scaled in the same way as QLE. In addition to the scaling relationships defined in \eqref{scaling}, an additional reference is needed for the current, we choose $\bar{j} = \tfrac{\mu k_{B}T_{0}\bar{n}}{Le}$ and set $j' = j/\bar{j}$. Using these conventions, the scaled QDD model has the following form (again, omitting the primes on the dimensionless variables): with $\nabla=d/dx$,
\begin{equation}
\label{qdd}
\left\{
    \begin{split}
        &\partial_t n + \nabla (n\nabla(\Ascr+\Vscr+\Vscr^{\rm{ext}}))=0 \\[2mm]
        &\alpha^2\Delta  \Vscr = n, \quad \Vscr(0)=\Vscr(1)=0, \\
        & n = n[e^{-\Ham_{\Ascr(t)}}]=\sum_{p\in\mathbb{N}} e^{-\lambda_p[\Ascr(t)]}\abs{\phi_p[\Ascr(t)]}^2,
    \end{split}\right.
\end{equation}
where $\{\lambda_p[\Ascr(t)], \phi_p[\Ascr(t)]\}_{p \in \Nm}$ are the spectral elements of the Hamiltonian $\Ham_{\Ascr(t)} = -\beta^2\Delta +\Ascr(t) = \Ham_0 + \Ascr(t)$. As with QLE, the Hamiltonian is equipped with Neumann boundary conditions. Finally, insulating boundary conditions are specified for the electrochemical potential $\Ascr+\Vscr+\Vscr^{\rm{ext}}$, i.e. \[\frac{d}{dx}(\Ascr+\Vscr+\Vscr^{\rm{ext}})\vert_{x=0,1} = 0.\] 
With such conditions, the total number of particles is preserved in the domain and there is no current at the boundary, as for QLE. The relationship with the solution $\varrho(t)$ to QLE is that $\varrho(t) \simeq \exp(-\Ham_{\Ascr(t)})$ as $\eps \to 0$.

Maybe counterintuitively,  QDD is probably best seen as an evolution equation on the chemical potential $\Ascr$ and the Poisson potential $\Vscr$ rather than on the density $n$. The mathematical analysis of \fref{qdd} is quite difficult, and an existence result in a one-dimensional periodic domain is obtained in \cite{pinaudpoincare}.
%\tr{As shown in \cite{QDDentropic}, \cite{frenchdocument}, the numerical method for the quantum drift diffusion model reduces to a problem of minimizing the free energy functional at each time step.  To obtain $\Ascr$ and $\Vscr$ at each step, the free energy 
%\[ \Fscr(\varrho) = \Tr(\varrho\log\varrho-\varrho)+\Tr(\varrho\Ham_0)+\Tr(\Vscr^{\text{ext}}\varrho)+\frac{1}{2}\lpnorm{\nabla \Vscr}{2}^2,\]
%is minimized under the constraint that the local density of the minimizer, $n_\rho$, is equal to a given local density, $n_\varrho$. NOT CONVINCED}

\section{Numerical method for QLE} \label{numQLE}
We introduce in this section the numerical scheme for  QLE. We start with the time discretization, and prove the convergence of a semi-discrete Strang splitting scheme to the solution to  QLE. We then define the spatial discretization in a second step, and detail the resolution of the transport and collision parts.
\subsection{Time discretization: Strang splitting}
We first consider a semi-discrete model by discretizing the time variable. As already mentioned, the main difficulty in the resolution of QLE is the calculation of the nonlinear term $\varrho_e[\varrho]$ in the collision part. The problem is considerably simplified by using a splitting approach: writing
\begin{equation*}
    i\varepsilon\partial_t \varrho = \mathcal{L}(\varrho) + i\mathcal{Q}(\varrho) :=  \frac{1}{\sqrt{2}\beta}[\Ham,\varrho]+\frac{i}{\varepsilon}(\rhoeq{\varrho}-\varrho),
\end{equation*}
we define two subproblems by splitting the operator on the right-hand-side into a transport part, $\mathcal{L}(\varrho)$, and a collision part, $\mathcal{Q}(\varrho)$. The collision subproblem is given by
\begin{equation}
\label{strang1}
    \eps \partial_t \varrho_1 =  \calQ(\varrho_1),\quad \varrho_1(t=0)=\varrho_1^{(0)},
  \end{equation}
  and the transport subproblem by 
\begin{equation}
\label{strang2}
    i\varepsilon\partial_t \varrho_2 = \calL(\varrho_2), \quad \varrho_2(t=0)=\varrho_2^{(0)}.
\end{equation}
Note that both problems are nonlinear since $\Ham$ involves the Poisson potential, and we have actually $\Ham\equiv \Ham[\varrho(t)]$. The latter is not difficult to handle compared to $\varrho_e[\varrho]$, and this is why it is included in the Hamiltonian part. %The decomposition $\calL(\varrho)+i Q(\varrho)$ can be seen in spirit as a decomposition into linear and nonlinear parts.
  
The crucial observation here is that \fref{strang1} preserves the local density (we write $n[\varrho]$ for $n_\varrho$ when it is more convenient): indeed, by construction of the equilibrium $\varrho_e[\varrho_1]$, we have $n[\varrho_e[\varrho_1]]=n[\varrho_1]$, and as a consequence, by linearity of the trace,  
\begin{equation*}
 \eps\partial_t n[\varrho_1] =   n[\calQ(\varrho_1)]= \frac{1}{\varepsilon}(n[\varrho_{e}[\varrho_1]]-n[\varrho_1]) = 0.
\end{equation*}
Hence, $n[\varrho_1(t)]=n[\varrho_1(0)]=n[\varrho_1^{(0)}]$, and the collision subproblem then becomes the linear equation
%\begin{equation} \label{lincol}
$$ 
   \partial_{t}\varrho = \frac{1}{\varepsilon^2}(\varrho_{e}[\varrho_{1}^{(0)}] -\varrho_{1}), \quad \varrho_{1}(t=0) = \varrho_{1}^{(0)}.
$$ 
 %\end{equation}
  We explain in Section \ref{secnumcol} how this problem is solved numerically.
% We begin by first looking at the model with time discretized and the spatial variable is continuous. The time-discrete solution to the quantum Liouville-BGK equation is obtained via a Strang operator splitting method. When devising a scheme for the QLE, the nonlinear collision term provides the most difficulty. The splitting scheme is set up into two subproblems, one with the nonlinear collision term and one with the free Liouville equation. Thanks to the local density conservation property of the collision operator, the nonlinear BGK-collision operator becomes perfectly linear which eases the numerical difficulties when solving the QLE. To begin, we rewrite the QLE as follows
% \begin{equation*}
%     i\varepsilon\partial_t \varrho = \mathcal{L}(\varrho) + i\mathcal{Q}(\varrho) =  \frac{1}{\sqrt{2}\beta}[\Ham,\varrho]+\frac{i}{\varepsilon}(\rhoeq{\varrho}-\varrho),
% \end{equation*}
% and define two subproblems by splitting the operators on the right-hand-side into a linear part, $\mathcal{L}(\varrho)$, and a nonlinear part, $\mathcal{Q}(\varrho)$. The linear part is the free Liouville subproblem given by
% \begin{equation}
% \label{strang2}
%     i\varepsilon\partial_t \varrho_2 = \frac{1}{\sqrt{2}\beta}[\Ham,\varrho_2] = \calL(\varrho), \quad \varrho_2(t=0)=\varrho_2^0,
% \end{equation}
% and the nonlinear part is the collision subproblem given by
% \begin{equation}
% \label{strang1}
%     \partial_t \varrho_1 = \frac{1}{\varepsilon^2}(\varrho_e[\varrho_1]-\varrho_1) = i\calQ(\varrho),\quad \varrho_1(t=0)=\varrho_1^0.
%   \end{equation}
  
We now express the Strang splitting scheme. The solution to each subproblem \fref{strang1} and \fref{strang2} can formally be represented in terms of an evolution operator, i.e.
\begin{equation*}
    \varrho_1(t) = W(t)\varrho_1^{(0)}, \qquad \text{and} \qquad \varrho_2(t) = U(t)\varrho_2^{(0)}.
\end{equation*} %Note that the solution to the transport subproblem is given by $\varrho_{2}(t)=U(t)\varrho_{2}^{0} = e^{-it\Ham}\, \varrho_{2}^{0}\, e^{it\Ham}$. % Thus, given an initial density operator, the Strang splitting solution at time $t$, denoted  is given by
% \begin{equation}
% \label{strangsplitscheme}
%     \varrho_{s}(t) = e^{\calL t/2}e^{\calQ t}e^{\calL t/2}\varrho^{0}.
% \end{equation}
For $h>0$, let $t_{k}=k\dt$ for $k=0,1,2,\cdots$. For a given initial condition $\varrho^0$, the semi-discrete Strang solution at time $t+t_{k-1}$, denoted $\varrho_{s}(t+t_{k-1})$,  is then obtained from the solution at $t_{k-1}$ by, for $k \geq 1$,
\[ \varrho_{s}(t+t_{k-1}) = U(t/2)W(t)U(t/2)\varrho_{s}^{k-1}, \qquad t \in [0,\dt], \qquad \varrho_s^0=\varrho^0, \]
% Let $\varrho_{s}^0 = \varrho(t=0)$ be the given initial condition.
with $\varrho_s^{k-1}=\varrho_{s}(t_{k-1})$. Thus, the Strang solution at time $t_k$ is given by
\[ \varrho_{s}^{k} = U(\dt/2)\underbrace{W(\dt)U(\dt) \cdots W(\dt)U(\dt)}_{k-1\text{ times}}W(\dt)U(\dt/2)\varrho^0.\]
We show in the next section that this scheme is well-defined and converges to the continuous solution as $h \to 0$. The important point to check is that the collision subproblem \fref{strang1} can indeed be solved at each time step. This amounts to verify that the solution $\varrho_s^k$ satisfies adequate conditions at each $k$.

\subsection{Convergence analysis} \label{convanal}
%We show in this section that the splitting scheme is convergent. 
We do not prove optimal estimates in the time step parameter $h$ since the optimal regularity of the map $\varrho \mapsto \varrho_{e}[\varrho]$ is still an open problem. It is known so far that the map has H\"older regularity $1/8$ in the space of Hilbert-Schmidt operators, which is enough for our purpose of showing convergence of the scheme. Moreover, we are not interested here in the asymptotic properties of the scheme as $\eps \to 0$, and will therefore set $\eps=1$ in the proof to simplify notation. The constant $C$ in the estimate of our convergence Theorem \ref{conv} further then depends on $\eps$ and grows as $\eps$ decreases to 0.

We first recall the existence result of \cite{MP-JDE} for the quantum Liouville-BGK equation. Note that the result therein is stated for the free Schr\"odinger operator, that is without any potentials. We will therefore set the Poisson and the external potentials to zero in this section to be consistent with \cite{MP-JDE}. We believe though that the result of \cite{MP-JDE} can be directly adapted to include these potentials (and as a consequence so does our convergence result below), but this is beyond the scope of this work.

Before stating the result, we need to introduce a few functional spaces. The space $\calJ_1$ is the space of trace class operators on $L^2(0,1)$  with norm $\jnorm{\varrho}=\Tr (|\varrho|)$, where $\abs{\varrho} = \sqrt{\varrho^*\varrho}$ for $\varrho^*$ the adjoint of $\varrho$; and $\calJ_{2}$ is the space of Hilbert-Schmidt operators on $L^2(0,1)$ with norm $\hsnorm{\varrho} = (\Tr(\varrho^{*}\varrho))^{1/2}$. The space $\calH$ is defined as
 \[\mathcal{H} = \{\varrho \in \calJ_{1}, \text{ such that } \overline{\Ham_{0}\abs{\varrho}\Ham_{0}}\in\calJ_{1}\},\]
where $\overline{\Ham_{0}\abs{\varrho}\Ham_{0}}$ denotes the extension of the operator $\Ham_{0}\varrho \Ham_{0}$ to $L^2(0,1)$, it is a Banach space when equipped with the norm 
\[\hnorm{\varrho} = \Tr(\abs{\varrho})+\Tr(\overline{\Ham_{0}\varrho\Ham_{0}}).\]
In the same way, $\calE$ is the space  \[\mathcal{E} = \{\varrho \in \calJ_{1}, \text{ such that } \overline{\sqrt{\Ham_{0}}\abs{\varrho}\sqrt{\Ham_{0}}}\in\calJ_{1}\},\]
and is Banach when equipped with the norm
\[\hnorm{\varrho} = \Tr(\abs{\varrho})+\Tr(\overline{\sqrt{\Ham_{0}}\varrho \sqrt{\Ham_{0}}}).\]
We will drop the extension sign in the sequel for simplicity. The space $\calE_+$ is the space of nonnegative operators in $\calE$, and we recall that a density operator is a self-adjoint, trace class, nonnegative operator.
  The result of \cite{MP-JDE} is the following:

 \begin{theorem} \label{exist}Suppose that the initial density operator $\varrho^0$ is in $\calH$, is such that $\varrho^0=f(\Ham_0)+\delta \varrho$, $f(\Ham_0) \in \calE_+$, $\delta \varrho$ self-adjoint in $\calE$, and that there exists $\underline{n}>0$ such that
   $$
   n[f(\Ham_0)](x) \geq \underline{n}, \qquad \forall x \in [0,1], \qquad \textrm{and} \qquad \|\delta \varrho\|_{\calE} \leq \underline{n}/4.
   $$
   Then, for any $T>0$, the QLE equation admits a solution $\varrho$ in $C^0([0,T],\calH) \cap C^1([0,T],\calJ_1)$ satisfying the integral equation
   \begin{equation}
\label{qleintsln}
    \varrho(t) = e^{-t}U(t)\varrho^0 + \int_0^t e^{-(t-s)}U(t-s)\rhoeq{\varrho(s)}ds,
  \end{equation}
where $U$ is the solution operator to the free Liouville equation (with $\Vscr=\Vscr^{\rm{ext}}=0$) introduced in the previous section. Moreover, the density verifies
  %\be \label{lb}
$$
  n[\varrho(t)](x) \geq e^{-T} \underline{n}/2, \qquad \forall (t,x) \in [0,T]\times[0,1].
$$ 
 %\ee
  
\end{theorem}

Note that the above result only provides us with the existence of solutions. We will actually prove the uniqueness further, by comparing any solution to the integral equation \fref{qleintsln} to the unique density operator obtained by the splitting scheme. Theorem \ref{exist} is actually stated in \cite{MP-JDE} in the context of periodic boundary conditions, and holds for the Neumann boundary conditions considered here with minor modifications.

To obtain the integral representation of the splitting solution and compare it with the original solution, we use the fact that the solution to the collision subproblem \fref{strang1}, $\varrho_{1}(t) = W(t)\sigma$ is given by
 \begin{equation*}
   \varrho_{1}(t)=W(t)\sigma =e^{-t}\sigma + \int_0^t e^{-(t-s)}\varrho_e[\sigma]ds.
 \end{equation*}
 
 Given $T$ and $h \leq 1$ positive, we denote by $N_T$ the largest integer such that $N_T h \leq T$. Thus, denoting by $\varrho_{s}^{k}$ and $\varrho^k$ the splitting solution and a solution to the integral equation at time $t_k= k h$, respectively, we have, for $t \in [0,h]$,
\begin{equation} \label{intsplit}
    \varrho_{s}(t_k+t) = U(t/2)W(t)U(t/2)\varrho^{k}_s= e^{-t}U(t)\varrho^{k}_{s} + \int_0^t e^{-(t-u)} U(t/2)\rhoeq{U(t/2)\varrho^k_{s}}du,
\end{equation}
and
\begin{equation} \label{intex}
    \varrho(t_k+t) = e^{-t}U(t)\varrho^{k} + \int_0^t e^{-(t-s)} U(t-s)\rhoeq{\varrho(t_k+s)}ds.
\end{equation}
For $t\in[0,\dt]$, let $e_{k}(t) := \varrho(t_k+t)-\varrho_{s}(t_k+t)$, where again $\varrho$ is \textit{any} solution to the integral equation \fref{qleintsln}. Note that we have by definition $e_{k+1} = e_{k}(h)$. 

The result below, proved in Section \ref{proofsplitHbnd}, shows that the splitting solution is well-defined and bounded in $\calH$.

\begin{lemma} \label{splitHbnd}
   Under the conditions of Theorem \ref{exist} on $\varrho^0 \in \calH$, the splitting scheme admits a unique nonnegative solution in $\calH$ with the following bound
\begin{equation} \label{eq:unfHbnd}
    \hnorm{\varrho^{k}_{s}}\leq  e^{C k h} \hnorm{\varrho^0},  \qquad \forall k \geq 0,
  \end{equation}
  where $C$ is a constant independent of $k$ and $h$. 
Furthermore, the splitting scheme preserves the trace, i.e.
%\begin{equation} \label{eq:traceconserved} 
$$ 
   \jnorm{\varrho^{k}_{s}} = \jnorm{\varrho^0}, \qquad \forall k \geq 0,
$$ 
 %\end{equation}
  and the local density $n[U(\tau)\varrho_s^k]$ verifies
  $$
  n[U(\tau)\varrho_s^k] \geq e^{-T} \underline{n}/2, \qquad \forall \tau \geq 0.
  $$
\end{lemma}

The next lemma, proved in Section \ref{prooflocalerror}, provides us with a local error estimate.
 \begin{lemma} \label{localerror}
 %Suppose that the relationship between a given density operator and its corresponding equilibrium operator, $\varrho \mapsto \rhoeq{\varrho}$, is H{\"o}lder continuous in the  $\calJ_{2}-$norm with exponent $\gamma=1/8$. 
Under the conditions of Theorem \ref{exist}, the local error $e_k$ between a solution to \fref{intex} and the splitting solution satisfies, for each $k$ and all $t \in[0,h]$,
% \begin{equation}\label{eq:localerroreqn}
 $$
    \hsnorm{e_{k}(t)} \leq \hsnorm{e_{k}} + C(h^{\frac{1}{1-\gamma}}+h^{1+\gamma} + h^2),
$$
 %  \end{equation}
   where $\gamma=1/8$ and the constant $C$ is independent of $h$ and $k$.
 \end{lemma}

Iterating the local estimate of Lemma \ref{localerror}, we arrive at the following result.

\begin{theorem} \label{conv}
  Under the condition of Theorem \ref{exist} on $\varrho^0 \in \calH$, we have, for any solution $\varrho$ to the integral equation \fref{qleintsln},
  $$
  \hsnorm{\varrho^{N_{T}}-\varrho^{N_{T}}_{s}} \leq C(h^{\frac{\gamma}{1-\gamma}}+h^{\gamma} + h),
  $$
  where $\gamma=1/8$, $C$ is independent of $h$, and $\varrho_s$ is the splitting solution.
\end{theorem}
Indeed, according to Lemma \ref{localerror},

\begin{equation*}
    \hsnorm{\varrho^{N_{T}}-\varrho^{N_{T}}_{s}}\leq \hsnorm{e_{N_{T}-1}} + C\big(\dt^\frac{1}{1-\gamma}+\dt^{1+\gamma}+ \dt^2\big), \end{equation*}
and iterating yields the desired estimate
\begin{equation*}
\begin{split}
    \hsnorm{\varrho^{N_{T}}-\varrho^{N_{T}}_{s}}& \leq C(h^{\frac{1}{1-\gamma}}+h^{1+\gamma} + h^2) N_T \leq C(h^{\frac{\gamma}{1-\gamma}}+h^{\gamma} + h).
\end{split}
\end{equation*}

At that point, we have therefore obtained that the unique splitting solution is close to any solution to the QLE for small $h$. Note that the error estimate of Theorem \ref{conv} is by no means optimal, as mentioned at the beginning of the section. If the map $\varrho \mapsto \varrho_e[\varrho]$ is Lipschitz (which we believe holds but cannot prove yet), then we expect as usual with Strang splitting to find a global order of convergence of two. This fact was verified numerically. %Note that when discretizing the spatial variable, it is possible to show that the now finite-dimensional map $\varrho \mapsto \varrho_e[\varrho]$ is Lipschitz with a Lipschitz constant that depends on the number of discretization points, and therefore that the fully discrete scheme is of order two in time for a given spatial discretization. We verify this fact numerically in Section blablabla.

\paragraph{Uniqueness for the continuous equation.}
A by-product of Theorem \ref{conv} is the uniqueness of solutions of \fref{qleintsln}. Fix indeed some $t>0$ and $t>h>0$, and write $t=N_t h+r_h$, with $N_t \in \Nm$ and $r_h \in [0,1)$. Consider then two possible solutions to \fref{qleintsln}, denoted $\varrho_1$ and $\varrho_2$. The associated splitting solution $\varrho_s$ is unique and verifies, according to Theorem \ref{conv} and Lemma \ref{localerror},
$$
\| \varrho_j(t)-\varrho_s(t)\|_{\calJ_2} =o(h), \qquad j=1,2.
$$
Hence, by the triangle inequality,
$$
\| \varrho_1(t)-\varrho_2(t)\|_{\calJ_2}=o(h),
$$
and since both $t$ and $h$ are arbitrary, this means that $\varrho_1=\varrho_2$ for all $t$. Uniqueness for nonlinear PDEs is often obtained under a Lipschitz condition on the nonlinearity, which, as mentioned, has not been established here. Uniqueness for our problem is a consequence of three factors: (i) the fact that the minimizer $\varrho_e[\varrho]$ is unique for a given $\varrho$, yielding a unique splitting solution, (ii) the equation for the collision part of the splitting scheme becomes linear, and (iii) the H\"older regularity of the map $\varrho \mapsto \varrho_e[\varrho]$.

Since the exact solution is now unique, we then conclude from Theorem \ref{conv} that the splitting solution converges to the unique solution to \fref{qleintsln}.

We  now turn to the spatial discretization of QLE.

\subsection{Spatial discretization}

Since we will compare the solutions to the QLE and QDD equations, we use the same spatial discretization for both, and adopt the one proposed for QDD in \cite{QDD-SIAM}. We discretize the (nondimensionalized) spatial domain $[0,1]$ with $N+2$ points $x_{p} = p \dx$ for $p = 0,1,\dots, N+1$ and $\dx = 1/(N+1)$. For a smooth function $\varphi$, integrating $\Ham \varphi$ over the interval $[x_{p-1/2},x_{p+1/2}]$ for $1 \leq p \leq N$ yields
\begin{align*}
\frac{1}{\dx} \int_{x_{p-1/2}}^{x_{p+1/2}}&\Ham \varphi (x) dx\\
&=-\frac{\beta^2}{\dx}\left( \varphi'(x_{p+\frac{1}{2}})-\varphi'(x_{p-\frac{1}{2}})\right)- \frac{1}{\dx} \int_{x_{p-1/2}}^{x_{p+1/2}}(\Vscr+\Vscr^{\text{ext}})(x)\varphi (x) dx\\
&=-\frac{\beta^2}{\dx^2}\left( \varphi(x_{p+1})-2 \varphi(x_{p})+\varphi(x_{p-1})\right)-(\Vscr+\Vscr^{\text{ext}})(x_p)\varphi (x_p)+O(\dx^2).
\end{align*}
Above, we used the midpoint rule for the integral. Note that we make sure when setting the discretization that the exterior potential $\Vscr^{\text{ext}}$ is smooth in each interval $(x_{p-1/2},x_{p+1/2})$. Since $\Vscr^{\text{ext}}$ typically has jumps, the discretization is chosen such that the jumps occur at some of the midpoints $x_{p+1/2}$ and not in $(x_{p-1/2},x_{p+1/2})$.

As in \cite{QDD-SIAM}, we adopt a first order discretization of the Neumann boundary conditions, resulting in $\varphi(x_0)=\varphi(x_1)$ and $\varphi(x_N)=\varphi(x_{N+1})$, and in the discrete $N \times N$ Neumann Laplace operator

\[ \Delta_{\text{Neu}}= \frac{1}{\dx^2}
    \begin{bmatrix}
    -1 & 1 & 0 & \cdots & 0 \\
    1 & -2 & 1 & \cdots & 0 \\
    \vdots & \ddots & \ddots &  & \vdots \\
    0 & \cdots & 1 &-2 & 1 \\
    0 & \cdots & \cdots & 1 & -1
    \end{bmatrix}.
    \]
The overall order of the spatial scheme is therefore one. The discrete Dirichlet Laplace operator used for the calculation of the Poisson potential has the same expression as $\Delta_{\text{Neu}}$, with the first and last entries on the diagonal replaced by $-2/\dx^2$.

Integrals are approximated in the same manner as in \cite{QDD-SIAM} as follows:
\bee \int_0^1 \varphi(x) dx &=& \int_0^{x_{1/2}} \varphi(x) dx+\int_{x_{N+1/2}}^{x_{N+1}} \varphi(x) dx+\sum_{p=1}^N\int_{x_{p-1/2}}^{x_{p+1/2}} \varphi(x) dx\\
&=&\int_0^{x_{1/2}} \varphi(x) dx+\int_{x_{N+1/2}}^{x_{N+1}} \varphi(x) dx+\dx \sum_{p=1}^{N}\varphi(x_p)+O((\dx)^2)\\
&=&\dx \sum_{p=1}^{N}\varphi(x_p)+O(\dx).
\eee
The boundary integrals above are discarded since the Neumann boundary conditions are accounted for at first-order only. The discrete inner product on $\Cm^N$ is then, for two vectors $u,v$,
$$\langle u,v \rangle = \dx \sum_{p=1}^{N}\bar{u}_pv_p.$$ %This allows for the definition of the discrete Euclidean 2-norm, $\norm{u}_{2} = (u,u)^{1/2}$. 
With the discrete Laplacian operator, we can now define the discrete analogs of the different Hamiltonian operators:
\begin{equation*}
    H_0 = -\beta^2\Delta_{\text{Neu}}, \quad H = H_{0} -\text{diag}(V+V^{\text{ext}}), \quad H_{A} = H_{0} + \text{diag}(A),
\end{equation*}
where $\text{diag}(w)$ for $w\in \mathbb{R}^{N}$ denotes the diagonal matrix with vector $w$ on the diagonal (we will just write $H_0+W$ for $H_0+\text{diag}(W)$ to simplify) and $V^{\text{ext}}$, $V$, and $A$ are the discrete counterparts to the exterior, Poisson, and chemical potentials, respectively. The Poisson equation becomes
\be \label{disP}
\alpha^2 \Delta_{\rm{Dir}} V=n_\varrho, 
\ee
with
$$
n_\varrho=\sum_{p=1}^N \rho_p |\phi_p|^2 \in \Rm^N,
$$
for $\{\rho_p, \phi_p\}_{1\leq p \leq N}$ the eigenvalues and eigenvectors of the positive matrix $\varrho$. All discrete eigenvectors are normalized such that $\langle \phi_p,\phi_p \rangle=1$.

Since the convergence of the semi-discrete splitting has already been established, and this is the most difficult part, it is a standard matter to prove that the fully discretized scheme is convergent. We omit the details.

\subsubsection{The collision subproblem} \label{secnumcol}
In this section, we detail the resolution of the collision subproblem \eqref{strang1}. We recall it has the following form: 
\begin{equation} \label{coleq}
    \partial_{t}\varrho_1 = \frac{1}{\varepsilon^2}(\varrho_{e}[\varrho_{1}^{(0)}] -\varrho_{1}), \quad \varrho_{1}(t=0) = \varrho_{1}^{(0)},
\end{equation}
where  $\rhoeq{\varrho}$ denotes the minimizer of the discrete free energy 
\[F(\sigma) = \TR(\sigma \log\sigma -\sigma) + \TR(H_0\sigma),\]
over nonnegative matrices $\sigma$ such that $n_\sigma=n_\varrho$. Note that compared to \fref{FE}, it is enough to consider $\TR(H_0\varrho)$ in $F$ instead of $\TR(H\sigma)$ since $\TR(H\sigma) = \TR(H_0\sigma)-\langle V+V^{ext},n_\sigma \rangle$, and the second term is fixed as $n_\sigma=n_\varrho$. The equilibrium is actually not calculated by minimizing $F$ under constraints, but rather by exploiting the form of the minimizer \fref{quantM}, and by the unconstrained minimization of the nondimensional discrete equivalent of the functional $\mathscr{J}(A)$ defined in \fref{defJ}, that is
\be \label{disJ} J(A) = \sum_{p=1}^{N}e^{-\lambda_{p}[A]} + \langle A,n_\varrho\rangle,\ee
where $\{\lambda_p[A]\}_{1 \leq p \leq N}$ is the set of eigenvalues of $H_A$. It is proved in \cite{QDD-SIAM} that the functional $J(A)$ is strictly convex and admits therefore a unique minimizer. The minimization procedure for $J(A)$ is described in detail in the next section.
% Due to this convention, the equilibrium operator is given by 
% \[\rhoeq{\varrho_{1}^0} = e^{-(H_0+A[\varrho_{1}^{0}])},\]
% where the chemical potential $A[\varrho_{1}^{0}]$ is chosen such that $n_{\rhoeq{\varrho_1^0}} = n_{\varrho_1^0}$. The chemical potential is obtained by performing an unconstrained minimization on the functional $J(A)$, which is equivalent to the Lagrangian for the constrained free energy minimization. 
% \begin{equation}
%     \label{JA}
%   J(A) = \sum_{p=1}^{N}e^{-\lambda_{p}[A]} + (A,n_{\varrho_1^0}),
% \end{equation}
% where $(\cdot,\cdot)$ denotes the discrete inner product between the given local density and the chemical potential $A$. The minimization procedure for $J(A)$ is described in detail in the next section.

Once $A[\varrho_{1}^{(0)}]$ is obtained by the minimization procedure, and therefore the equilibrium operator $\varrho_{e}[\varrho_{1}^{(0)}]=\exp(- (H_0+A[\varrho_{1}^{(0)}]))$ is known, the now linear collision problem \fref{coleq} is reduced to a set of coupled ODEs that describe the evolution of the operator $\varrho_1$. The solution is easily found to be
$$
\varrho_1(t)=\expt{t} \varrho_1^{(0)} +(1-\expt{t})\varrho_{e}[\varrho_{1}^{(0)}].
$$
From a practical viewpoint, the matrices $\varrho_1^{(0)}$ and $\varrho_{e}[\varrho_{1}^{(0)}]$ are defined on different basis of $\Rm^N$. We then express them both in the canonical basis to form $\varrho_1(t)$, and diagonalize the resulting matrix to store the spectral elements of $\varrho_1(t)$.

% % To do so, we must consider the spectral representation of both the initial density operator and its associated equilibrium operator. 
% Let indeed $\varrho_1^{(0)}$ be given by
% \[\varrho_1^{(0)} = \sum_{p=1}^N \sigma_{p}\ket{\phi_p}\bra{\phi_p},\] 
% with its associated equilibrium operator
% \[\varrho_{e}[\varrho_1^{(0)}] = \sum_{p=1}^N e^{-\lambda_p[A]} \ket{\psi_p[A]}\bra{\psi_{p}[A]},\]
% with $\{\lambda_{p}[A], \psi_{p}[A]\}_{1 \leq p \leq N}$ the spectral elements of the discrete Hamiltonian $H_{A} = H_0+A[\varrho_{1}^{(0)}]$. We then simply project $\varrho_{1}$ onto the basis $\{\psi_{p}[A]\}_{1 \leq p \leq N}$ leading to
% \[ \varrho_1(t) = \sum_{p=1}^N\mu_p(t)\ket{\psi_p[A]}\bra{\psi_p[A]},\] 
% and
% \begin{equation*}
%     \mu'_{p}(t) = \frac{1}{\varepsilon^2}(e^{-\lambda_p[A]}-\mu_p(t)), \quad \mu_{p}(t=0) = \bra{\psi_p[A]}\varrho_{1}^{(0)} \ket{\psi_p[A]}, \quad \text{for } p = 1, \dots, N,
% \end{equation*}
% and the solution is given by 
% \[ \mu_{p}(t) = e^{-\lambda_p[A]} -\big(e^{-\lambda_p[A]} - \mu_p(t=0)\big)\expt{t},\]
% for $p=1,\dots,N$. 

We describe in the next section the minimization of the functional $J(A)$.
% The nondimensional discrete equivalent of the functional $\mathscr{J}(A)$ defined in \fref{defJ} is then given by
% \[ J(A) = \sum_{p=1}^{N}e^{-\lambda_{p}[A]} + \langle A,n_\varrho\rangle,\]
% where $\varrho$ is the splitting solution to the 
% The discrete equilibrium operator can be represented in terms of the eigenvalues and eigenvectors of the discrete Hamiltonian $H_{A}$,
% \[\rhoeq{\varrho} = \sum_{p=1}^{N} e^{-\lambda_{p}[A]}\ket{\psi_{p}[A]}\bra{\psi_{p}[A]},\]
% where the chemical potential $A$ in the Hamiltonian is obtained by minimizing the functional $J(A)$. 
\subsubsection{Minimization procedure}
 We use the Polak-Ribi{\`e}re variant of the nonlinear conjugate gradient algorithm to minimize $J(A)$. For a given local density $n$ (replace $n_\varrho$ by $n$ in \fref{disJ}), the unique minimizer $A_\star$ is such that $n_{e^{-H_{A_\star}}}=n$. We start with an initial guess $A^{(0)}$, and must find an initial search direction $s^{(0)}$ and step length $b^{(0)}$ to initialize the algorithm. We set $s^{(0)}= -\nabla_A J(A^{(0)}) \in \Rm^N$, with, see e.g. \cite{QDD-SIAM}, 
\begin{equation*}
    \nabla_{A}J(A) = n-n_{e^{-H_A}} = n-\sum_{p=1}^{N}e^{-\lambda_{p}[A]}\abs{\psi_{p}[A]}^2 ,
\end{equation*}
for $\{\lambda_{p}[A], \psi_{p}[A]\}_{1 \leq p \leq N}$ the spectral elements of  $H_{A} = H_0+A$. We find the step length via a line search 
 $$
 b^{(0)} = \underset{b \in \Rm} {\rm{argmin}}\; J(A^{(0)} + b\, s^{(0)}).$$
 We will see further that it is possible to obtain a very good initial guess for the line search, and, as consequence, a simple method avoiding the calculation of the Hessian such as the secant method proves to be efficient. Once $b^{(0)}$ is found, we update the chemical potential as $A^{(1)} = A^{(0)} + b^{(0)} s^{(0)}$. \\

% The line search is performed via Newton's method, with the Hessian of $J(A)$ given explicitly by, see \cite{QDD-SIAM},
% $$
% \left(\nabla_A^2 J(A)\right)_{ij} = \sum_{p=1}^{N}\sum_{q=1}^{N} \omega_{pq}\Psi_{pq}^i(\Psi_{pq}^j)^*, \qquad i,j=1, \cdots,N,$$
% where
% \begin{equation*}
% \begin{split}
%     \omega_{pq}& = \begin{cases} &-\frac{e^{-\lambda_p[A]}-e^{-\lambda_q[A]}}{\lambda_{p}[A]-\lambda_{q}[A]} \quad \text{if  } p\neq q \\ &e^{-\lambda_{p}[A]}  \quad \text{ if  } p=q,
% \end{cases}, \qquad \qquad 
% \Psi^i_{pq} = \psi^i_{p}[A]\psi^i_{q}[A],
% \end{split}
% \end{equation*}
% for $\psi^i_{p}[A]$ the $i$-th component of the vector  $\psi_{p}[A] \in \Cm^N$.

 The nonlinear conjugate gradient algorithm is then as follows:\\

 % While  $\|dA^{(k)}\|_{\ell_1} > \text{tolerance}$:
 While  $\|A^{k}-A^{k-1}\|_{\ell_2}/\|A^{k}\|_{\ell_2} > \text{tolerance}$:
 \begin{itemize}
     \item Compute the steepest descent direction, $dA^{(k)} = -\nabla_A J(A^{(k)})$,
     \item Compute $c = \max\{0,c_{PR}\}$ where $c_{PR} = \frac{\langle dA^{(k)}, dA^{(k)}-dA^{(k-1)}\rangle}{\langle dA^{(k-1)},dA^{(k-1)} \rangle}$,
     \item Update the search direction $s^{(k)} = dA^{(k)}+c\,s^{(k-1)}$,
     \item Perform line search $b_k = \underset{b \in \Rm} {\rm{argmin}}\;J(A^{(k)} + b\,s^{(k)})$,
     \item Update chemical potential $A^{(k+1)} = A^{(k)}+b^{(k)}s^{(k)}$.
 \end{itemize}

We explain in Section \ref{simu} how the algorithm can be accelerated by exploiting some particular regimes of parameters, in particular one where $\beta$ is small.
 
We consider next the resolution of the transport part \fref{strang2} of the splitting scheme.

\subsubsection{The transport subproblem}
We recall that the spatially discrete version of \eqref{strang2} is 
\begin{equation*}
    i\varepsilon\partial_t \varrho_2 = \frac{1}{\sqrt{2}\beta}[H,\varrho_2], \qquad \varrho_2(t=0)=\varrho_2^{(0)}=\sum_{p=1}^{N}\gamma_p\ket{v_p}\bra{v_p},
\end{equation*}
where $H = -\beta^2\Delta_\textrm{Neu} -V^{\rm{ext}}-V$, for $V\equiv V[\varrho_2]$ the Poisson potential and $\{\gamma_{p}, v_{p}\}_{1 \leq p \leq N}$ the spectral elements of $\varrho_2^{(0)}$. The solution to the above system is given by
\[ \varrho_{2}(t) = \sum_{p=1}^{N}\gamma_{p}\ket{v_p(t)}\bra{v_p(t)}, \] 
where $v_p(t)$ is the solution to the nonlinear Schr\"odinger equation 
\bee
i \partial_t v_p(t) &=& \frac{1}{\sqrt{2}\beta\varepsilon}H v_p(t) = \frac{1}{\sqrt{2}\beta \varepsilon}\left(-\beta^2\Delta_{\text{Neu}} -V^\text{ext}\right)v_p(t) - \frac{1}{\sqrt{2}\beta \varepsilon}V(t)v_p(t)\\
&=:& H_Lv_p(t) + B(t)v_p(t),\eee
with initial condition $v_p(t=0)=v_p$. The above equation is nonlinear because of $V(t)\equiv V[\varrho_2(t)]$, and becomes linear when using Strang splitting for the time discretization as for the collision term. The approximate Strang solution $v_p^{(1)}$ for the above nonlinear Schr\"odinger equation at time $t=h$ is then given by
 \[v_{p}^{(1)} = e^{-i\dt H_L/2} S(h) e^{-i\dt H_L/2}v_p,\]
where $w_p(t):=e^{-itH_{L}}w_p^{(0)} $ is the solution to the linear Schr{\"o}dinger equation 
 \be \label{eqw}i\partial_t w_{p}(t) = H_Lw_p(t) = \frac{1}{\sqrt{2}\beta\varepsilon}(-\beta^2\Delta_{\text{Neu}} - V^\text{ext})w_p(t), \quad \text{for } p = 1,\dots, N,\ee
with $w_p(0)=w_p^{(0)}$, and where  the second subproblem $z_p(t):=S(t)z_p^{(0)} $ is reduced to the following set of ODEs
 \be \label{eqz}
i\partial_t z_p(t) = -\frac{1}{\sqrt{2}\beta\varepsilon} V(t) z_p(t), \quad \text{for } p = 1, \dots, N.
\ee
In terms of the density operator $\varrho_2$, the splitting scheme yields at $t=h$  the approximate solution
$$
\varrho_2^{(1)}=e^{-i\dt H_L/2} S(h) e^{-i\dt H_L/2} \varrho_2^{(0)}e^{i\dt H_L/2} S^*(h) e^{i\dt H_L/2}.
$$
The key point is that the equation \fref{eqz} on $z_p$ preserves the absolute value of $z_p$ since $V$ is real-valued. This means that for a density operator $\sigma$, we have for all $t \geq 0$,
$$
n[S(t) \sigma S^*(t)]=n[\sigma].
$$
Therefore, the Poisson potential $V(t)$ in \fref{eqz} is actually linear, equal to
$$V(t)=V[S(t) e^{-i\dt H_L/2} \varrho_2^{(0)}e^{i\dt H_L/2} S^*(t)]=V[e^{-i\dt H_L/2} \varrho_2^{(0)}e^{i\dt H_L/2}]=V(0),$$
and is the solution to $\alpha^2 \Delta_{\rm{Dir}} V=n[e^{-i\dt H_L/2} \varrho_2^{(0)}e^{i\dt H_L/2}]$.
 %where $V$ is the Poisson potential that solves $-\Delta_{dir} V = n[\varrho_2]$. 
 % We obtain an approximate solution $v_{p}^{k+1}$ at time step $t_{k+1}$ from the solution $v_{p}^{k}$ at time step $t_k$ by solving the linear system given by the Crank-Nicolson method, 

The solution to \fref{eqw} at time $t=h$ is obtained with the standard Crank-Nicolson scheme
 \be \label{crank} \left(i\mathbb{I}-\frac{\dt}{2} H_L\right)w^{(1)}_{p} = \left(i\mathbb{I}+\frac{\dt}{2} H_L\right)w_{p}^{(0)},\ee
 where $\mathbb{I}$ denotes the $N$ by $N$ identity matrix. % The second sub-problem is reduced to the following ODEs,
 % \[i\partial_t v_p = \frac{1}{\sqrt{2}\beta\varepsilon} Vv_p, \quad \text{for } p = 1, \dots, N,\]
 % where $V$ is the Poisson potential that solves $-\Delta_{dir} V = n[\varrho_2]$. Once the Poisson potential is obtained by solving the linear system 
 %  \[-\alpha^2\Delta^{D}_{\text{Dir}}V = n[\varrho_2^0],\]
 %  the approximate solution at the next time step can be obtained from the ODEs.
%   Then the approximate solution at the next time step is then given by an exponential \[u_p^{n+1} = e^{-\frac{i\Delta t}{\sqrt{2}\beta \varepsilon}V}v_{p}^{n}, \quad \text{for } p = 1, \dots, N.\]

\section{Numerical method for QDD} \label{numQDD}
We modify here the method introduced in \cite{QDD-SIAM}. In \cite{QDD-SIAM}, both the chemical potential $A$ and the Poisson potential $V$ are treated implicitly. At a given time, the potential $V$ is a minimizer of a given functional, and $A+V$ of another one. These two functionals are then combined in \cite{QDD-SIAM} in a somewhat arbitrary manner to form a unique functional to minimize. While from a purely theorical viewpoint this poses no problem, this creates unnecessary difficulties in the minimization of this latter functional. A simple way to improve efficiency is to treat $V$ explicitly, and to keep an implicit scheme for $A$ resulting in the minimization of a functional acting only on $A$.

The spatial grid is identical to that of QLE. The  fully discrete scheme adapted from \cite{QDD-SIAM} is first-order both in time and space, and reads
% \begin{equation*}
%     \begin{split}
%         &\frac{n^{k+1}_{i}-n^{k}_{i}}{\dt}+\frac{1}{2\dx^2}(n^{k}_{i+1}+n^{k}_{i})\big((A_{i+1}^{k+1}+W_{i+1}^{k+1}) - (A_{i}^{k+1}+W_{i}^{k+1})\big) \\&\quad \quad \quad \quad \qquad- \frac{1}{2\dx^2}(n^{k}_{i}+n^{k}_{i-1})\big((A_{i}^{k+1}+W_{i}^{k+1}) - (A_{i-1}^{k+1}+W_{i-1}^{k+1})\big)=0
%     \end{split}
% \end{equation*}
% where $i = 1, \dots, N-1$. This equation is coupled to the Poisson equation, which utilizes the discrete Dirichlet Laplacian,
% \[ -\alpha^2 \Delta^{D}_{\text{Dir}} V^{k+1} = n^{k+1},\]
% where $n^{k+1} = \sum_{i=1}^{N-1} e^{-\lambda_i[A^{k+1}]}\abs{\varphi_{i}[A^{k+1}]}^2$. 
% As noted in \cite{QDDentropic}, this finite volume scheme for QDD can be compactly written with the help of matrices that are used to decompose the Neumann Laplacian into a product of two matrices, namely $\Delta^{D}_{\text{Neu}} = \widetilde{D_{1}}D_{2} = \widetilde{D_{2}}D_{1}$. Where the $(N-1)$ by $(N-1)$ matrices are given by
\begin{equation} \label{qdddis}
  \left\{
    \begin{split}
        &\frac{n^{k+1}-n^{k}}{\dt} + \frac{1}{2}\tD{-}(n^{k}D_{+}(A^{k+1}+W^{k})) + \frac{1}{2}\tD{+}(n^{k}D_{-}(A^{k+1}+W^{k})) = 0 \\[3mm] 
        &\alpha^2 \Delta_{\rm{Dir}}V^{k} = n^{k}, \qquad W^{k}=V^k+V^{\rm{ext}}\\
        &n^{k+1} = \sum_{p=1}^{N} e^{-\lambda_p[A^{k+1}]}\abs{\psi_{p}[A^{k+1}]}^2,
    \end{split} \right.
\end{equation}
where $\{\lambda_p[A^k],\psi_p[A^k]\}_{1\leq p \leq N}$ are the eigenvalues and eigenvectors of $H_0+A^k$, and the $N\times N$ matrices  $D_+, D_-, \widetilde{D}_+, \widetilde{D}_-$ are given by 
$$
    D_+= \frac{1}{\dx}
    \begin{bmatrix}
    -1 & 1 & 0 & \cdots & 0 \\
    0 & -1 & 1 & \cdots & 0 \\
    \vdots & \ddots & \ddots &  & \vdots \\
    0 & \cdots & \cdots &-1 & 1 \\
    0 & \cdots & \cdots & 0 & 0
  \end{bmatrix}
 \qquad 
D_-= \frac{1}{\dx}
    \begin{bmatrix}
    0 & 0 & 0 & \cdots & 0 \\
    -1 & 1 & 0 & \cdots & 0 \\
    \vdots & \ddots & \ddots &  & \vdots \\
    0 & \cdots & \cdots & 1 & 0 \\
    0 & \cdots & \cdots & -1 & 1
    \end{bmatrix}
    $$

    and

    $$
    \tD{+}= \frac{1}{\dx}
    \begin{bmatrix}
    -1 & 1 & 0 & \cdots & 0 \\
    0 & -1 & 1 & \cdots & 0 \\
    \vdots & \ddots & \ddots &  & \vdots \\
    0 & \cdots & \cdots &-1 & 1 \\
    0 & \cdots & \cdots & 0 & -1
    \end{bmatrix}
  \qquad
     \tD{-}= \frac{1}{\dx}
    \begin{bmatrix}
    1 & 0 & 0 & \cdots & 0 \\
    -1 & 1 & 0 & \cdots & 0 \\
    \vdots & \ddots & \ddots &  & \vdots \\
    0 & \cdots & \cdots & 1 & 0 \\
    0 & \cdots & \cdots & -1 & 1
    \end{bmatrix}.
    $$

    The Neumann boundary conditions are accounted for in the definition of the above matrices, and the notation $UV$ in \fref{qdddis} for two vectors $U,V$ in $\Rm^N$ denotes the term-by-term product, i.e. $(UV)_i=U_iV_i$. Adapting \cite{QDD-SIAM}, given $n^k$ (and therefore $V^{k}$), the solution $A^{k+1}$ to the implicit problem \fref{qdddis} is obtained as the unique minimizer of the strictly convex functional
\begin{equation*}
    \begin{split}
J_{\textrm{QDD}}(A) 
=& \frac{\dt\dx}{4} \sum_{i=1}^{N} n_i^k(D_+(A+V^k+V^{\rm{ext}}))_i^2 + \frac{\dt\dx}{4} \sum_{i=1}^{N} n_i^k(D_-(A+V^k+V^{\rm{ext}}))_i^2 \\
& \quad \quad+ \sum_{i=1}^{N} e^{-\lambda_i[A]}+ \dx\sum_{i=1}^{N} n_i^kA_i.
\end{split}
\end{equation*}

% \noindent The QDD minimization scheme is given by the following steps:
% \begin{itemize}
%     \item Obtain an initial guess for chemical potential $A_0$ and Poisson potential $V_0$ given an initial density operator $\varrho^0$.
%     \item Update chemical potential $A^{k+1}$ and Poisson potential $V^{k+1}$ by minimizing $G(A,V)$ given the local density $n^{k}$ at time step $t_{k}$. 
%     \item Update the local density, $n^{k+1} = \sum_{p=1}^{N-1} e^{-\lambda_{p}[A^{k+1}]}\abs{\psi_{p}}^{2}$ where $(\lambda_p, \psi_p)$ are the spectral elements associated to $H = H_0 + A^{k+1}$.
% \end{itemize}
The minimization of the functional $J_{\textrm{QDD}}(A)$ is accomplished in the same manner as the collision step of the QLE, that is by using a nonlinear conjugate gradient method.

\section{Numerical results} \label{simu}
% In this section the numerical results for the QLE are presented. Since the QDD model is derived as a diffusive limit of the QLE, it is expected that in the long-time limit that the two models will agree and we will show that this holds true with the presented splitting scheme. In addition to comparing steady-state solutions of the QLE and the QDD, we compare the differences between the ballistic case (free Liouville equation) and the collisional case (the QLE). In the presented examples, the external potential takes on the roll of either a single or double barrier structure potential and we take the initial local density to be concentrated to the left of the potential barriers. 

\subsection{Complexity}

% Since the domain is one-dimensional and therefore no heavy calculations are required, we use Matlab for the simulations.
The resolution of the linear system \fref{crank} has to be repeated at each time step for the significant  modes in the density matrix; in the configurations we consider, there are between 50-100 modes used to build the density operator, and it turns out it is more effective to compute $\left(i\mathbb{I}-\frac{\dt}{2} H_L\right)^{-1}\left(i\mathbb{I}+\frac{\dt}{2} H_L\right)$ once and for all and then simply do the matrix vector multiplications. We use Matlab's backslash operator both for the inversion of $\left(i\mathbb{I}-\frac{\dt}{2} H_L\right)^{-1}$ and the resolution of the linear system \fref{disP} to obtain the Poisson potential $V$. The operator exploits the tridiagonal structure of the matrix for a cost of order $O(N)$.

% involved in the Crank-Nicolson scheme is done with Matlab's backslash operator which exploits the tridiagonal structure of the matrix for a cost of order $O(N)$. Note that it is actually more efficient to solve the linear system for each $p$ than computing the inverse matrix once and for all. The same method is applied to solve the linear system \fref{disP} to obtain the Poisson potential $V$.

The most expensive part of the simulation is the minimization of $J(A)$, which requires the (repeated) diagonalization of $H_0+A$. Since the matrices are fairly small in our simulations, say $500 \times 500$, it turns out it is actually faster to compute all eigenvalues with Matlab's \texttt{eig} function than using the function \texttt{eigs}, which computes only a small number of eigenvalues. Since \texttt{eig} is based on the QR method, and $H_0+A$ is already tridiagonal, the cost is $O(N^2)$ for each calculation of the eigenvalues.

\subsection{Initialization of the minimization algorithms} \label{init}
We need good initial guesses for best convergence of the nonlinear gradient algorithm for both $J(A)$ for QLE and $J_{\textrm{QDD}}(A)$ for QDD. They are obtained as follows. As mentioned at the end of Section \ref{subQLE}, the parameter $\beta$ is typically small in physically interesting regimes. It is then natural to exploit this fact to approximate $\varrho_{e,0}=\exp(-H_0+A^0)$ using semi-classical analysis. For the continuous problem, we show in Section \ref{semiC} in the Appendix that, for $x$ away from the boundaries,
\be \label{semiN}
    n[\exp(-\Ham_0+\Ascr)](x)= \frac{1}{\sqrt{4 \pi} \beta} e^{-\Ascr(x)}+o(1),
\ee
where $o(1)$ refers to a term that is small in appropriate sense when $\beta \ll 1$. As a consequence, we set as initial guess for the discrete problem $A_\textrm{guess}= - \log(\sqrt{4 \pi} \beta) n^0$. The latter provides a good approximation of the exact solution for $x$ away from the boundaries. 

At the time step $k$, we simply use the result $A^{k-1}$ of the previous step as initial guess.

\subsection{Acceleration of the nonlinear conjugate gradient}
\begin{figure}[h!]
    %\begin{center}
\centering
    \includegraphics[height=11cm, width=9.cm]{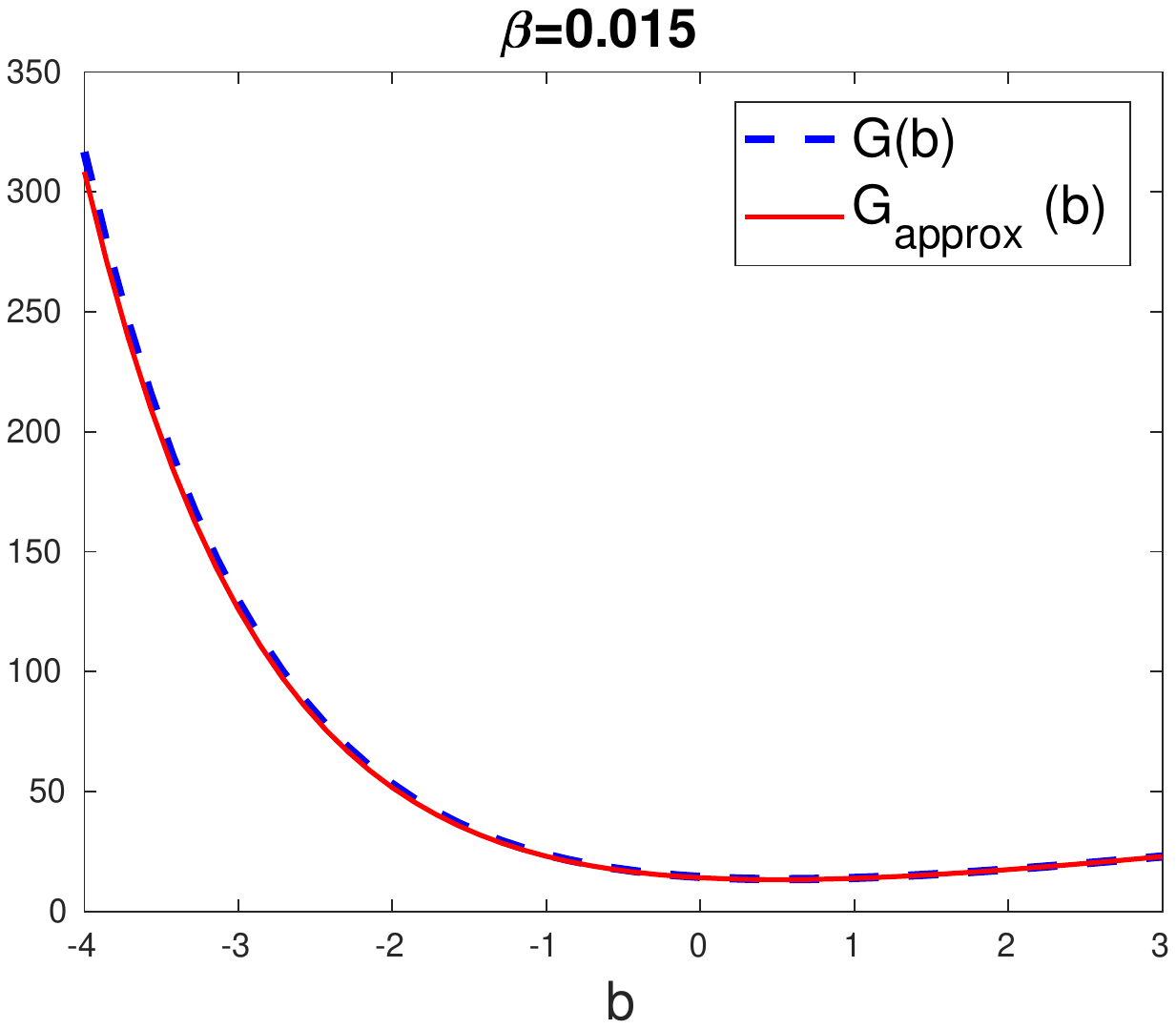} \hspace{-3cm}
    \includegraphics[height=11cm, width=9.cm]{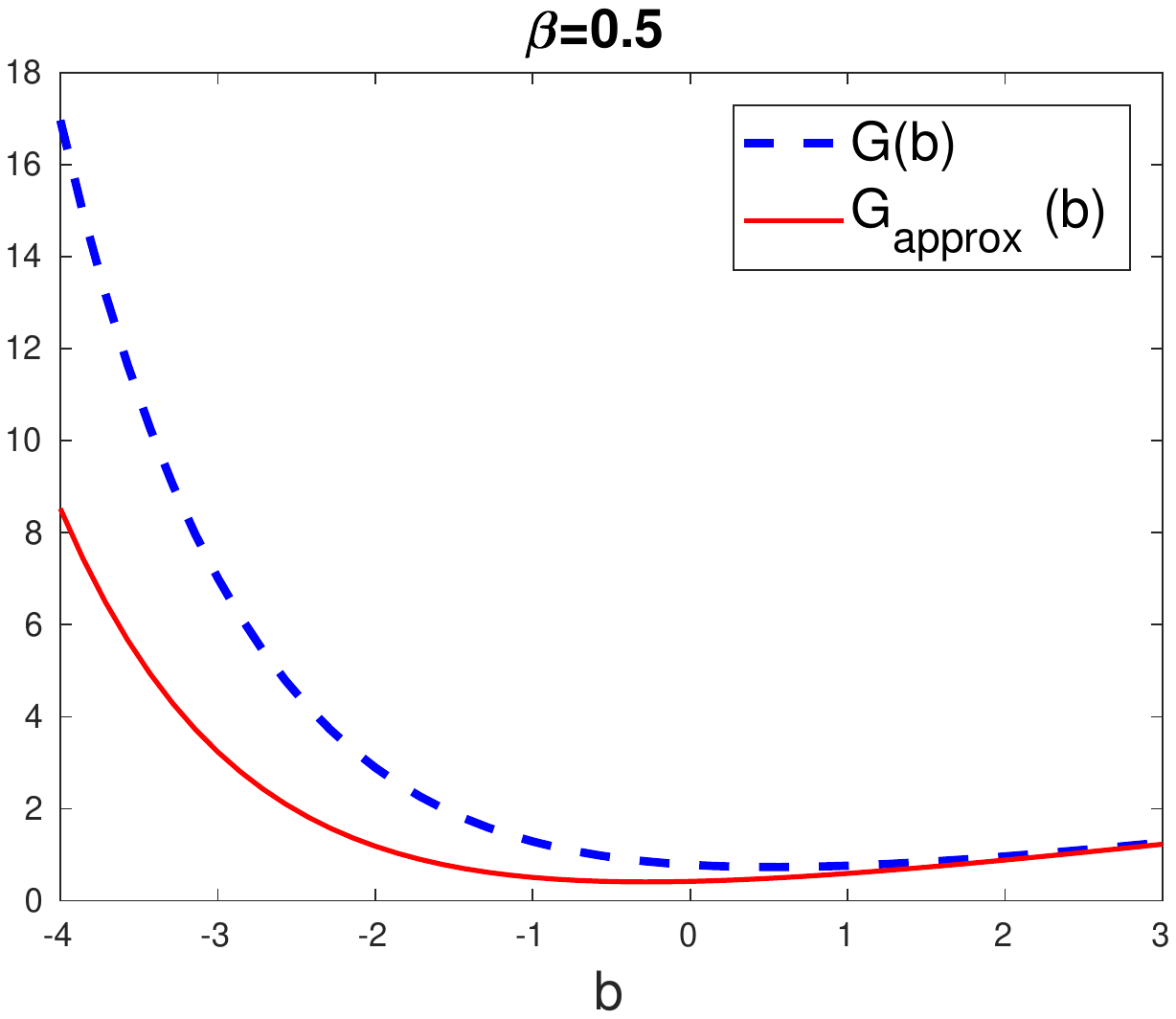}
%\end{center}    
\vspace{-3.3cm}
\caption{Comparison $G$ and $G_{\rm{approx}}$ for $A=\cos(4x)^3+x$, $s=1/(1+x^2)$, $n=n[\exp(-(H_0+A_1))]$ with $A_1=\cos((\cos(6x+1)))$.}
    \label{fig:G}
\end{figure}

Most of the computational time is spent in the diagonalization of the matrices $H_0+A$, and we explain here how to minimize the number of calls to the function \texttt{eig} in the minimization of $J(A)$ at each time step. As in the previous section, we exploit the fact that $\beta$ is small in our configuration of interest, and use \fref{semiN} to get an approximate expression of the functional $J(A)$. We then perform a line search with the approximate functional in order to get a good initial guess for the exact line search. For two vectors $A$ and $s$ given in $\Rm^N$, this approximate functional is shown in Section \ref{semiC} in the Appendix to be equal to, for $b \in \Rm$,
$$
G_{\rm{approx}}(b)=J_{\rm{approx}}(A+b s)=\frac{\Delta x}{\sqrt{4 \pi} \beta} \sum_{i=1}^N e^{-A_i+b s_i}+ b \langle A, s\rangle. 
$$
A straightforward Newton's method is used to find the minimizer of $G_{\rm{approx}}(b)$. While the function $G_{\rm{approx}}$ is not accurate for all values of $b$, it actually provides an excellent approximation of the minimizer of $G(b)=J(A+b s)$, even for values of $\beta$ up to 0.5, see figure \ref{fig:G}. Note that the behavior reported on the figure is not particular to the choice of $A$, $s$, and $n$, and holds for a large class of parameters.
\vspace{-3.3cm}

\subsection{Application: validation of QDD}

We compare in this section the solutions to QDD and to QLE for various values of $\eps$. We will see that the models agree when $\eps$ is sufficiently small, which is the regime of validity of QDD. We consider three situations: (i) in the first one, the initial condition is well-prepared in the sense that it is a quantum Maxwellian associated with a given Hamiltonian. This prevents the creation of initial layers as is customary in diffusion limits. We then switch at the initial time the potential in this Hamiltonian and observe how the system converges to a new equilibrium. (ii) The situation in the second case is slightly less favorable in the sense that the initial density operator is function of an Hamiltonian, but not a quantum Maxwellian. (iii) In the last scenario, we consider an ill-prepared initial condition that is a combination of wave packets; in this case, there is an initial layer and in order to minimize its effects and observe good agreement between QLE and QDD for earlier times, the parameter $\eps$ has to be decreased ($\eps=0.0025$ in the last case versus $\eps=0.01$ in the first two).

In all simulations, we set the tolerance for the nonlinear conjugate gradient and the associated line search to $10^{-7}$. The number of spatial discretization points is $N=400$. The parameter $\alpha$ is set to $\alpha=1$ for simplicity, which is of the order of magnitude of values found for semiconductor devices such as the resonant tunneling diode, for which $\alpha=1.7$ \cite{QDDscale}. As already mentioned, $\beta$ is small in interesting regimes, and we set for instance $\beta=0.015$. % in Section \ref{init}, the devices we are interested in are naturally in a semi-classical regime where $\beta$ is small, and we set for instance $\beta=0.015$.

In all density operators, we discard the modes associated with weights (i.e. eigenvalues) less than $10^{-7}$. This leaves approximately between 50 and 100 modes in the quantum Maxwellian for instance, and improves computational time. With $\beta=0.015$, and considering the quantum Maxwellian with the free Neumann Hamiltonian, we have about 50 modes with weights greater than $10^{-3}$, and about 30 others with weights between $10^{-3}$ and $10^{-7}$.

\begin{figure}[h!]
    %\begin{center}
\centering
    \includegraphics[height=5cm, width=6.cm]{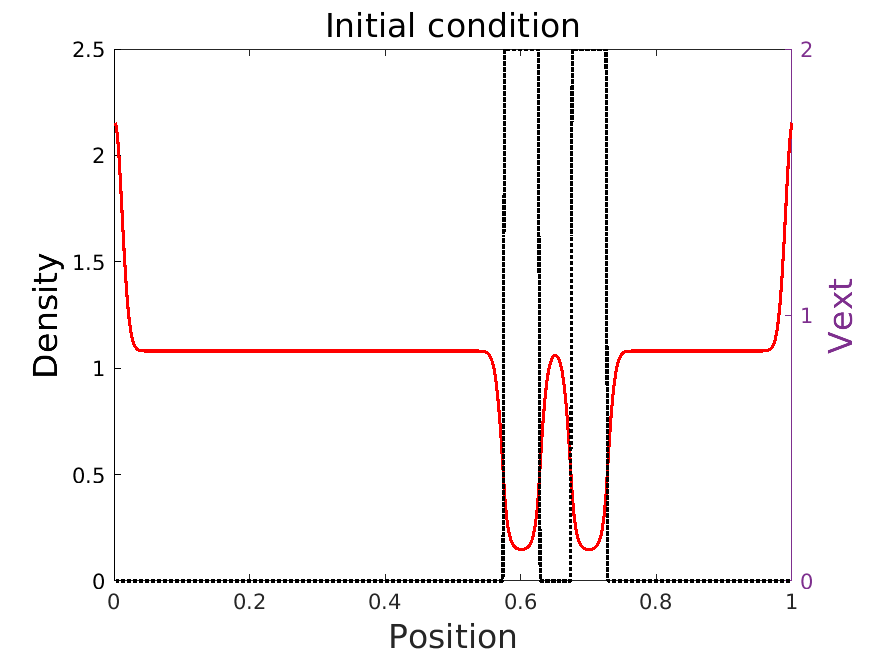} \hspace{0.5 cm}
    \includegraphics[height=5cm, width=6.cm]{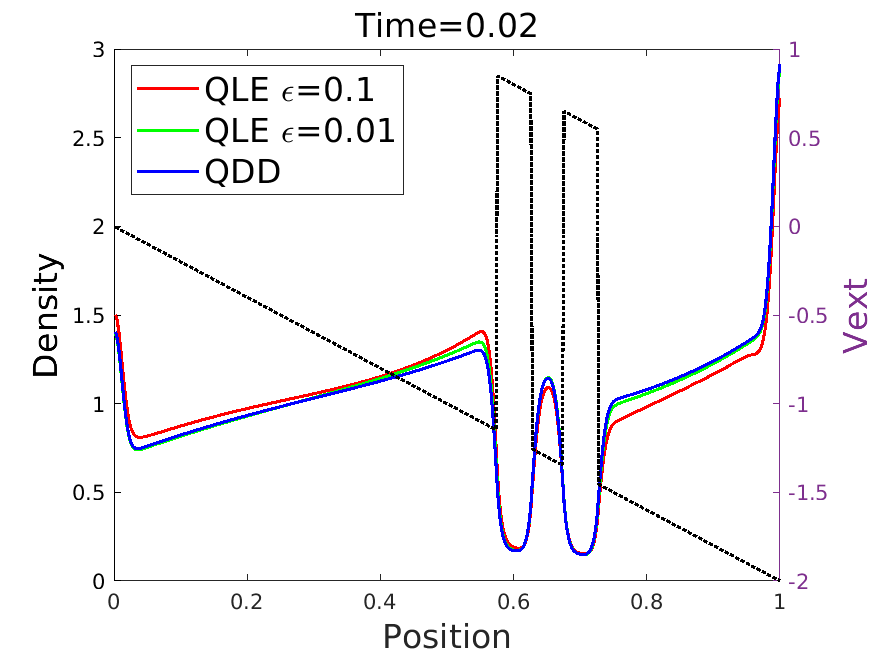}\vspace{0.5 cm}
    \includegraphics[height=5cm, width=6.cm]{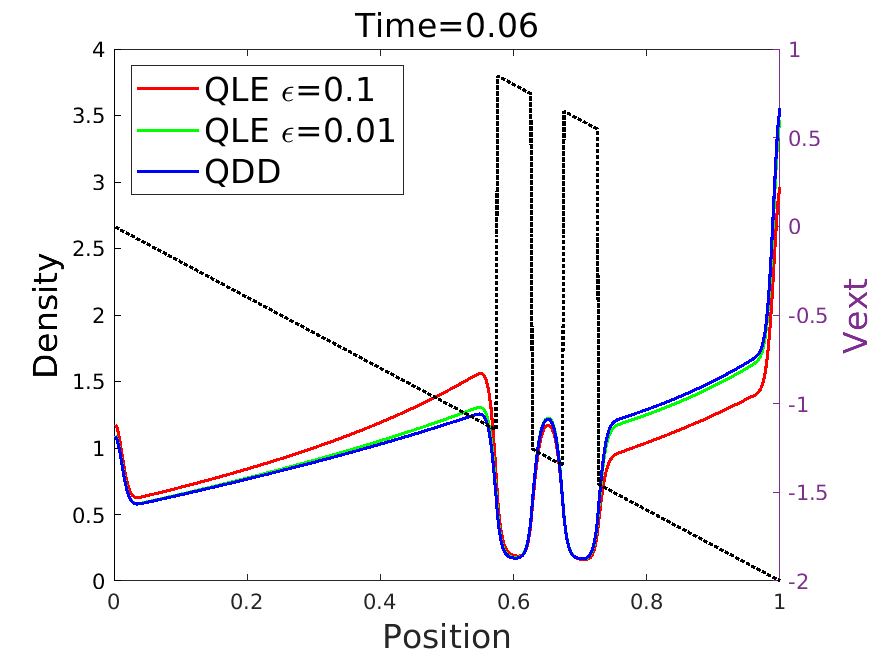}\hspace{0.5 cm} %\hspace{-3cm}
    \includegraphics[height=5cm, width=6.cm]{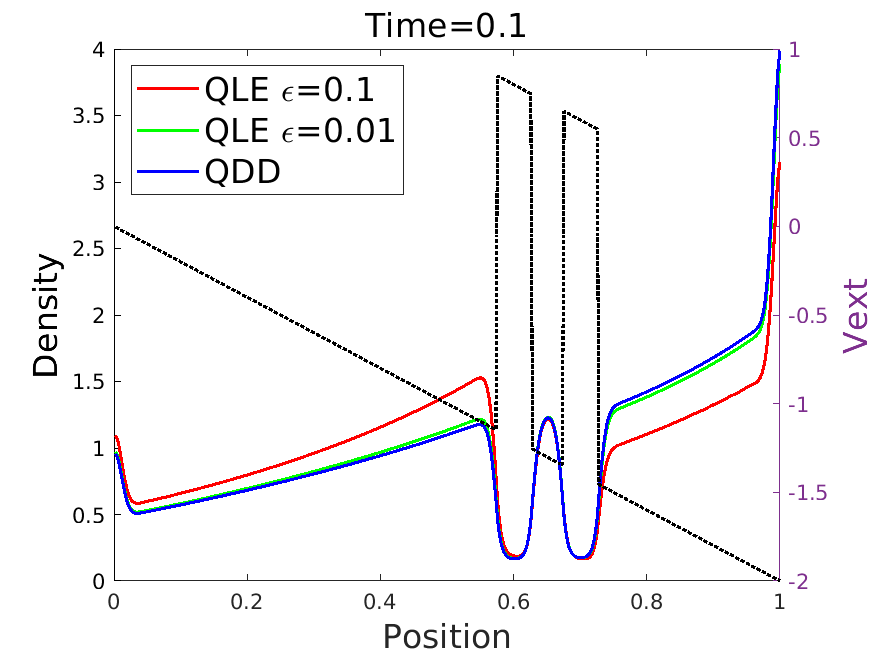}
%\end{center}    
%\vspace{-3.3cm}
\caption{Comparison QDD and QLE for the quantum Maxwellian initial condition. The initial double barrier potential is shifted at $t=0^-$ by $-2x$. The time $t=0.1$ is close to the equilibrium time for QDD. Observe the very good agreement when $\eps=0.01$ (error of order $2\%$).}
    \label{fig1}
\end{figure}

\paragraph{Quantum Maxwellian.} We set for initial condition

$$
\varrho_0=\frac{e^{-(-\beta^2 \Delta_{\textrm{Neu}}+V^{\textrm{ext,0}})}}{\Tr \big(e^{-(-\beta^2 \Delta_{\textrm{Neu}}+V^{\textrm{ext,0}})} \big)},
$$
where $V^{\textrm{ext,0}}$ is the double barrier potential shown in figure \ref{fig1}, top left panel (the width of the well and the barriers is 0.05, with height equal to 2). Such a potential is characteristics of the resonant tunneling diode, see \cite{QDDscale} and references therein. The density associated to $\varrho_0$ is depicted in the same panel. At time $t=0^-$, the potential $V^{\textrm{ext,0}}$ is switched to $V^{\textrm{ext}}=V^{\textrm{ext,0}}-2 x$, which is now the exterior potential used in the resolution of  QLE and QDD and promotes particle transport from left to right. It is depicted in the other panels of figure \ref{fig1}. The time stepsize is set to $h=10^{-4}$ for the calculations.

We then represent in figure \ref{fig1} the transition to the new equilibrium associated to $V^{\textrm{ext}}$, from time $t=0$ to $t=0.1$ (which is close to the time at which the equilibrium is reached by QDD). We observe a remarkable agreement between QDD and QLE with $\eps=0.01$, with an overall space-time relative $\ell_2$ error of about $2 \%$. When $\eps=0.1$, the diffusive regime is not valid and as a consequence QLE and QDD produce different densities.

\begin{figure}[h!]
    %\begin{center}
\centering
    \includegraphics[height=5cm, width=6.cm]{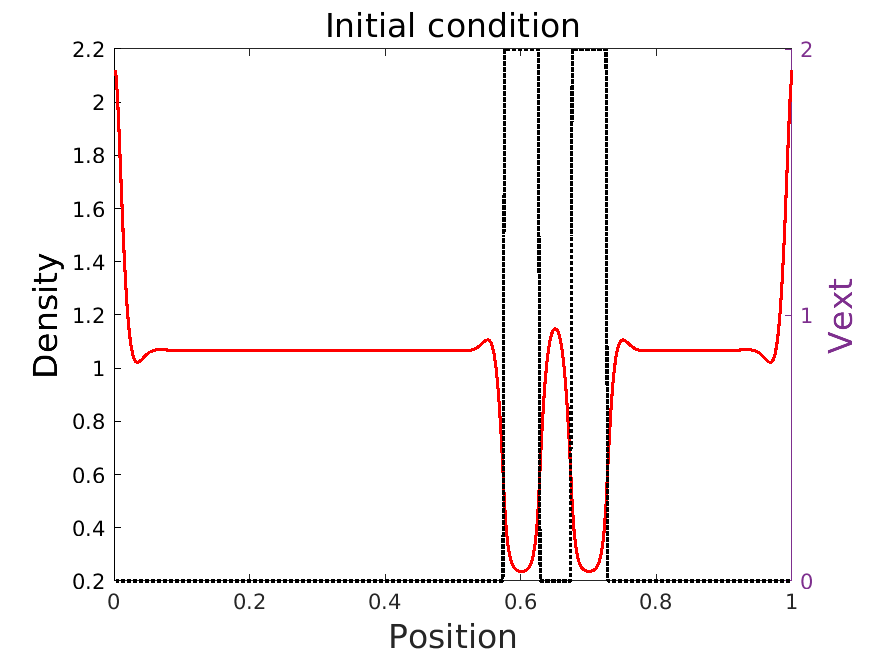} \hspace{0.5 cm}%\hspace{-3cm}
    \includegraphics[height=5cm, width=6.cm]{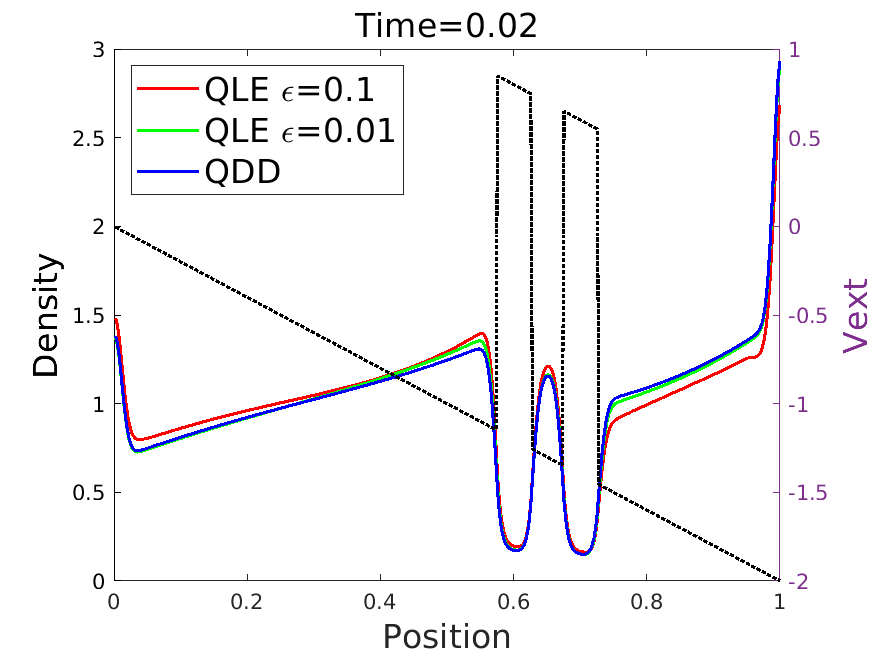}\vspace{0.5 cm}
    \includegraphics[height=5cm, width=6.cm]{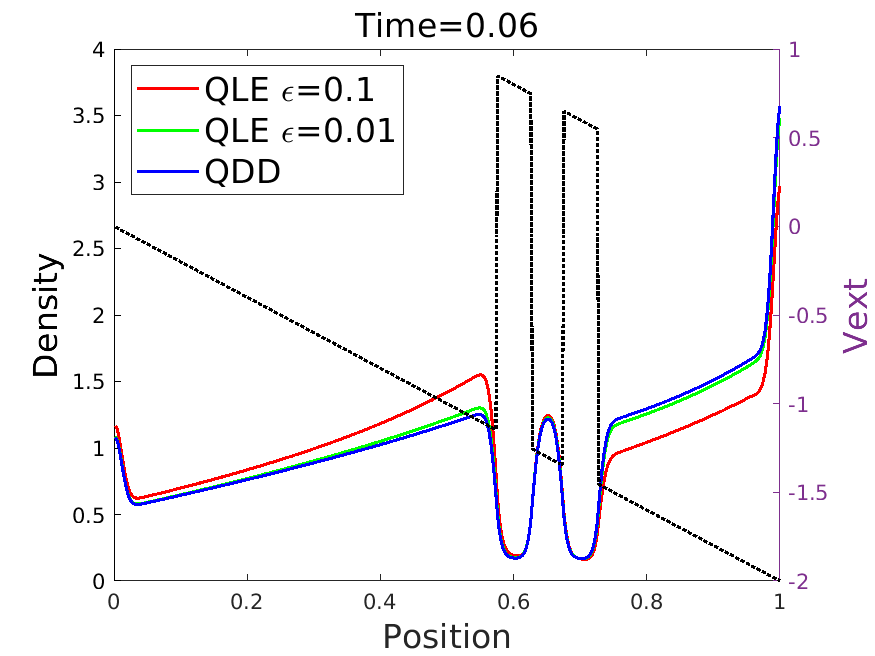} \hspace{0.5 cm}%\hspace{-3cm}
    \includegraphics[height=5cm, width=6.cm]{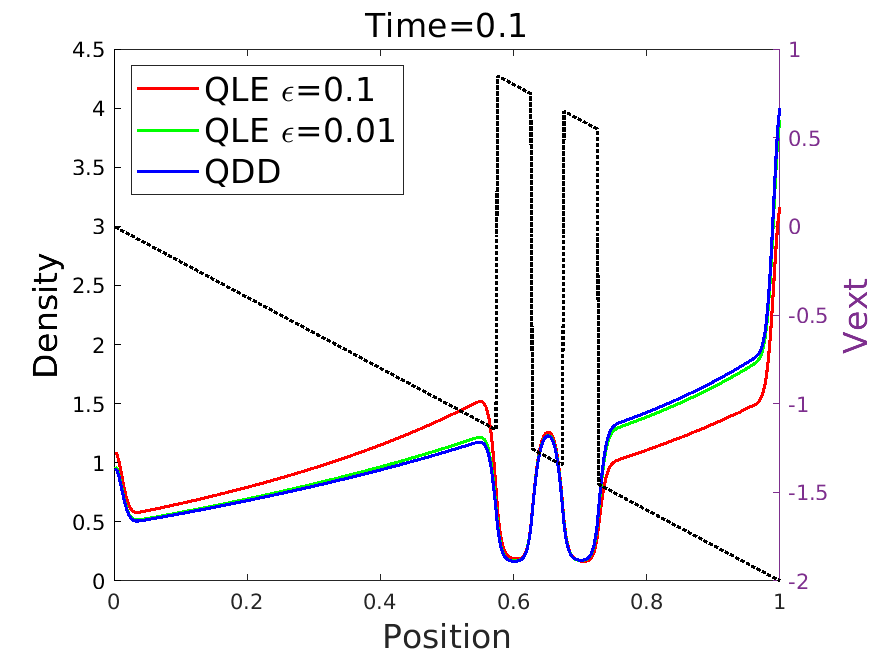}
%\end{center}    
%\vspace{-3.3cm}
\caption{Similar as figure \ref{fig1}, with now the initial condition given by a function of an Hamiltonian.}
    \label{fig2}
\end{figure}

\paragraph{Function of an Hamiltonian.} We set
$$
\varrho_0=\frac{f(-\beta^2 \Delta_{\textrm{Neu}}+V^{\textrm{ext,0}})}{\Tr \big(f(-\beta^2 \Delta_{\textrm{Neu}}+V^{\textrm{ext,0}}) \big)},
$$
with $f(x)=(1+x^2)^{-1}$ and the same parameters as in the previous paragraph. The situation is very similar as above with a very good agreement between QDD and QLE with $\eps=0.01$ and an error again of the order of $2 \%$. The densities are depicted in figure \ref{fig2}.

\begin{figure}[h!]
    %\begin{center}
\centering
    \includegraphics[height=5cm, width=6.cm]{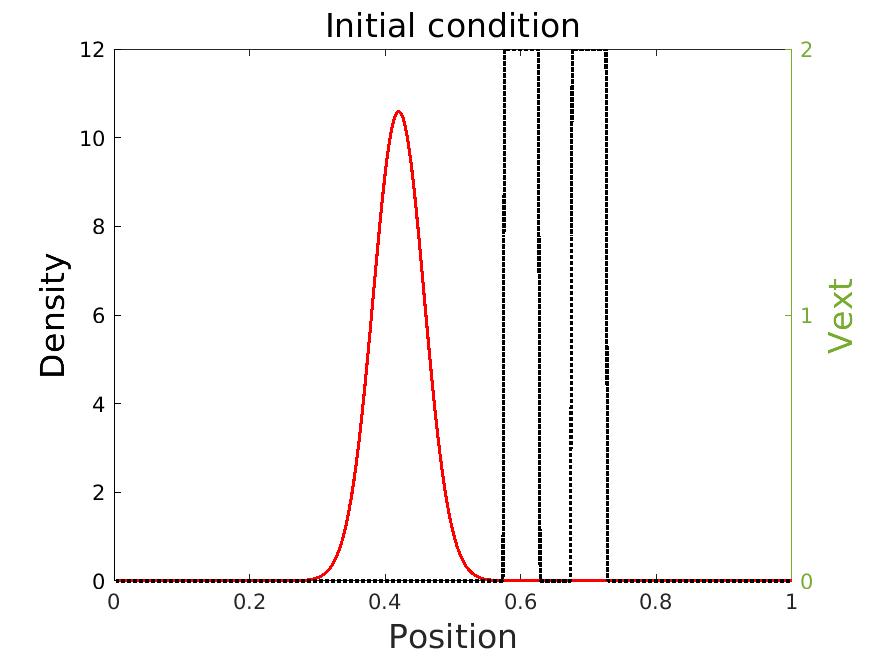} \hspace{0.5 cm}%\hspace{-3cm}
    \includegraphics[height=5cm, width=6.cm]{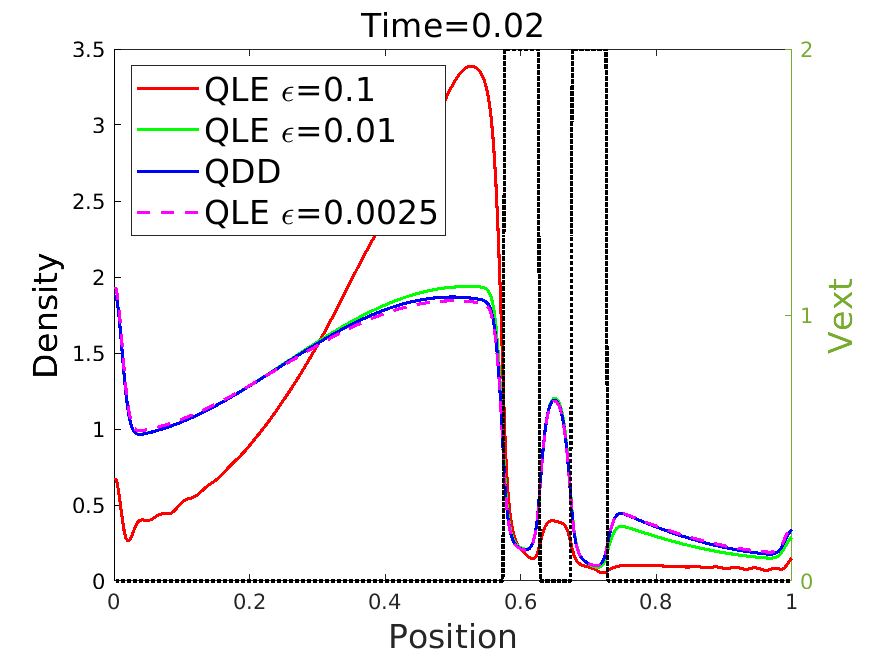}\vspace{0.5 cm}
    \includegraphics[height=5cm, width=6.cm]{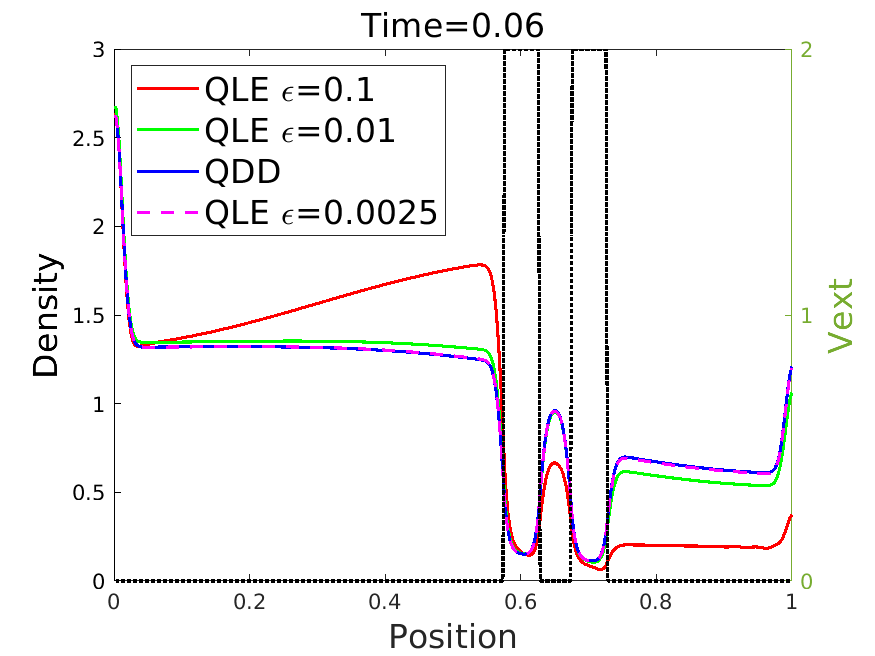} \hspace{0.5 cm}%\hspace{-3cm}
    \includegraphics[height=5cm, width=6.cm]{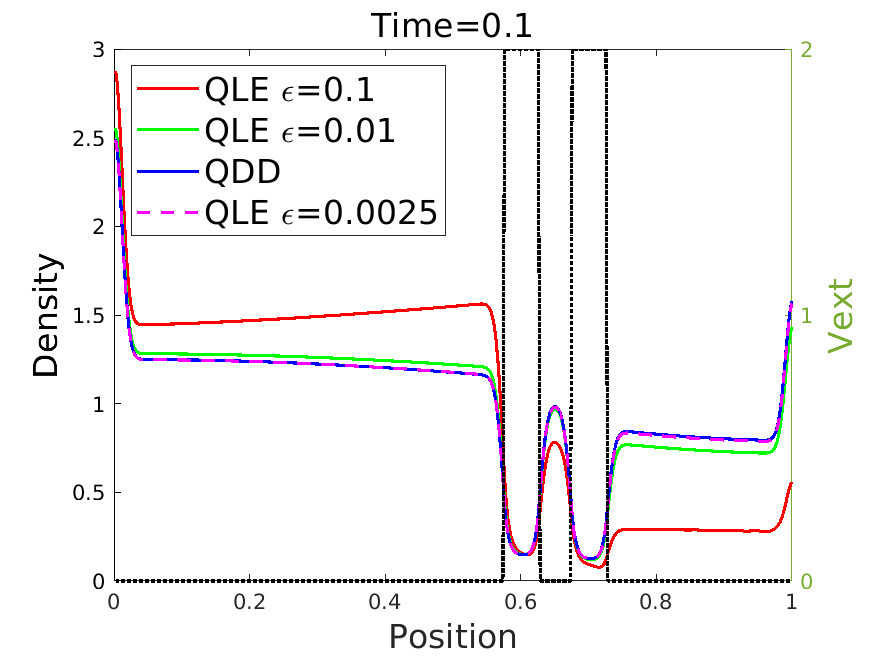}
%\end{center}    
%\vspace{-3.3cm}
\caption{Comparison QDD and QLE for the wave packets initial condition. The parameter $\eps$ has to be decreased to $\eps=0.0025$ to obtain a strong agreement (error of order $1\%$).}
    \label{fig3}
\end{figure}

\paragraph{Superposition of wave packets.} We set

$$\varrho^0=\frac{ \chi \gamma_0 \chi}{\Tr(\chi \gamma_0 \chi) }$$
where $\chi(x)$ is the function $\chi(x) = e^{-(x-x_0)^2/\sigma^2}+\eta$, and $\gamma_0$ is the density operator
$$
\gamma_0 = \sum_{p=1}^{5} e^{-\lambda_p}\ket{\psi_p}\bra{\psi_p},
$$
where $\lambda_p=(8 \pi \beta p)^2$ and $\psi_p(x) = e^{8 i \pi p x}$. The associated density is represented in the top left panel of figure \ref{fig3}, along with the (fixed this time) double barrier potential $V^{\textrm{ext}}=V^{\textrm{ext,0}}$ used in the calculations. In the localizing function $\chi$, we choose $x_0=0.42$ and $\sigma=0.075$. The parameter $\eta=5.10^{-3}$ acts as a regularization since the gaussian function is very small away from its center. Small densities create large chemical potentials $A$ which generate numerical instabilities, and we found that such an $\eta$ improves the convergence of the minimization algorithms.  

The simulations, represented in figure \ref{fig3}, show that $\eps=0.1$ is far too large to capture the diffusive regime. When $\eps=0.01$, the comparison improves with a space-time relative $\ell_2$ error of about $4 \%$. In order to observe a very good agreement, we decrease $\eps$ to $\eps=0.0025$. For the simulations with $\eps=0.0025$, the time stepsize is set at $5.10^{-6}$ to obtain sufficient accuracy. This substantially increases the numerical cost and makes the numerical method not effective for such small values of $\eps$. The relative error between QDD and QLE with $\eps=0.0025$ is now of order $1\%$. One would need to resort to asymptotic preserving schemes to capture the $\eps \ll 1$ regime at an affordable cost, see e.g. \cite{jinAP, lemouAP}.

% The initial data $\varrho^0$ for the QLE is a matrix, that is diagonalized to obtain the initial density $n^0$ for QDD. With $n^0$, one can then compute the initial Poisson potential $V^0$. It remains to find an initial guess for $A$. For this, let $\varrho_{e,0}$ be the minimizer associated with $n^0$. We look an approximation of $A^0$, the latter being such that $\varrho_{e,0}=\exp(-H_0+A^0)$. 
% We compare in this section the solutions to QLE and QDD. We choose the following parameters: we use $N_I=40$ modes in the initial condition $\gamma_0$, yielding a smallest wavelength $d_{\rm{min}}=2/N_I=1/20$. With $\beta=0.015$, the largest energy $\lambda_{N_I}=\beta^2 \pi^2 N_I^2 $ is 3.55, with a statistical weight of $e^{-\lambda_{N_I}}=0.0287$. Larger values of $\beta$ would decrease this weight and the resulting initial condition would essentially be low frequency. We set $\Delta x=d_{\rm{min}}/10$ for QLE. Regarding the time stepsize, we set $V^{\rm{ext}}=0$ in \fref{eqw}, so that the largest eigenvalue of $H_L$ contributing to the problem is of the order of $\lambda_{N_I}/(\sqrt{2} \beta \eps)$, resulting in time oscillations with period $T_{\rm{min}}=2 \pi \sqrt{2} \beta \eps / \lambda_{N_I}$. We then set $h=T_{\rm{min}}/10$.

%\subsection{Comparison between QLE and ballistic transport}

\section{Conclusion}

We have introduced in this work a time-splitting scheme for the resolution of the quantum Liouville-BGK equation. The splitting allows us, exploiting the local conservation of particles, to obtain a completely linear collision step. The minimization problem involved in the latter is solved by using the nonlinear conjugate gradient algorithm, and good initial guesses can be obtained by taking advantage of some small parameters. We applied our numerical method for comparing the solutions to the quantum Liouville-BGK equation and to the quantum drift-diffusion, and obtained excellent agreement in the regime of validity of the latter.

An important limitation of the method is the requirement that the time stepsize be small compared to the rescaled mean free path $\eps$ for good accuracy. We plan in the future on removing this restriction by designing an asymptotic preserving scheme in the spirit of \cite{jinAP, lemouAP}. This would allow us to capture the correct solution for arbitrarily small values of $\eps$ at a reasonable computational cost.
% \begin{table}
% \begin{center}
% \begin{tabular}{|c|c|c|c|c|}
%     \hline
%      $\Delta x$ &$\Delta t$ & $\alpha$ & $\beta$ & $\varepsilon$ \\
%      \hline
%      0.01 & 0.001 & 0.3334 & 0.116 &  0.05\\
%      \hline
% \end{tabular}
% \end{center}
% \caption{Constants used for the QLE and QDD steady state comparison.}
% \label{table:constsforcomp1}
% \end{table}

\begin{appendix}
  \section{Appendix}
  \subsection{Derivation of the 1D model} \label{deriv1D}
  
We derive in this section a 1D model as a simplification of a 3D model. 
We consider a 3D domain of the form $(0,L)\times \Omega$, where $\Omega$ is periodic as explained in the introduction, and we choose the 2-torus for simplicity. We then write
\be \label{tenL2}
L^2((0,L)\times \Omega) = L^2(0,L) \otimes L^2(\Omega),
\ee
in the sense that the two spaces are unitarily equivalent, and the 3D Hamiltonian is expressed as
$$
\mathscr{H}_{3D}= \mathscr{H}_x \otimes \un + \un \otimes \mathscr{H}_\perp,
$$
where (all physical constants are set to one),
$$
\mathscr{H}_x=-\frac{\partial^2}{\partial x^2}+W(x), \qquad x \in (0,L), \qquad \mathscr{H}_\perp=-\frac{\partial^2}{\partial y^2}-\frac{\partial^2}{\partial z^2}, \quad (y,z) \in \Omega,
$$
and $\un$ denotes, with an abuse of notation, the identity operator in both $L^2(0,L)$ and $L^2(\Omega)$. Above, $W$ is a given bounded potential. The operator $\mathscr{H}_x$ is equipped with the domain defined in \fref{domain}, and $\mathscr{H}_\perp$ with the domain consisting of $H^2(\Omega)$ periodic functions. The 3D Liouville-BGK equation is then
\begin{equation} \label{liou3D}
    i  \partial_t\varrho = [\mathscr{H}_{3D},\varrho] + i  (\varrho_e[\varrho]-\varrho), \qquad \varrho(t=0)=\varrho_0,
\end{equation}
where $\varrho_e[\varrho]$ is the unique minimizer of the 3D free energy
\begin{equation*}
    \Fscr_{3D}(\sigma) = \TR(\sigma \log \sigma) + \TR(\Ham_{3D} \sigma),
\end{equation*}
under the constraint that $n_\sigma=n_\varrho$. We set $\TR(\varrho_0)=1$, so that $\TR(\varrho(t))=1$ for all $t \geq 0$. Note that the traces above are taken w.r.t. $L^2((0,L)\times \Omega)$, and that we removed linear term $-\Tr(\sigma)$ in the entropy  since it is fixed to one by the constraint. Let
$$
\varrho_{\perp}=\frac{e^{-\mathscr{H}_\perp}}{\Tr_{\perp}(e^{-\mathscr{H}_\perp})},
$$
where $\TR_\perp$ denotes trace w.r.t. $L^2(\Omega)$. It is clear that $\varrho_\perp$ is the unique minimizer of the ``transverse'' free energy
$$
 \Fscr_{\perp}(\sigma)=S_\perp(\sigma)+\TR_\perp(\mathscr{H}_\perp \sigma),
$$
under the constraint that $\Tr_\perp (\sigma)=1$. The free energy $\Fscr_{\perp}(\sigma)$ is indeed, up to a constant term, equal to the relative entropy between and $\sigma$ and $\varrho_\perp$ which vanishes when $\sigma=\varrho_\perp$. Above, $S_\perp$ is the transverse entropy
$$
S_\perp(\sigma)=\TR_\perp (\sigma \log \sigma).
$$

We will show that if the initial condition $\varrho_0$ has the tensor form $\varrho_{0,x} \otimes \varrho_\perp$ ($\varrho_{0,x}$ acts on the space $L^2(0,L)$), namely that the initial state of the system is at equilibrium in the transverse plane, then the solution $\varrho(t)$ remains in a similar form and reads $\varrho(t)=\varrho_x(t) \otimes \varrho_\perp$. While it is direct to separate variables in $\Ham_{3D}$, it has to be proved that the minimizer $\varrho_e[\varrho(t)]$ for $\varrho(t)=\varrho_x(t) \otimes \varrho_\perp$ also admits a tensor form $\varrho_{e,x} \otimes \varrho_\perp$. This is a consequence of \fref{tenL2} and of the subaddivity of the von Neumann entropy $-\TR(\sigma \log \sigma)$. More precisely, we have the following lemma:

\begin{lemma} \label{subadd} Let $n(x,y,z)=n_0(x)/|\Omega|>0$, with $\|n_0\|_{L^1(0,L)}=1$. Then, the unique minimizer $\varrho^{3D}_\star$ of $\Fscr_{3D}(\sigma)$ with constraint $n_\sigma=n$ has the form
$$
\varrho_\star^{3D}=\varrho_\star[n_0] \otimes \varrho_{\perp},
$$
where $\varrho_\star[n_0]$ is the minimizer of the 1D problem with constraint $n_0(x)$.
\end{lemma}

The 1D problem mentioned in the lemma consists in minimizing the 1D free energy
$$
 \Fscr_{x}(\sigma)=S_x(\sigma)+\TR_x(\mathscr{H}_x \sigma),
$$
under the constraint that $n_\sigma=n_0$. Above, $S_x$ is the ``longitudinal'' entropy
$$
S_x(\sigma)=\TR_x (\sigma \log \sigma),
$$
where $\TR_x$ denotes trace w.r.t. $L^2(0,L)$. \bigskip

\begin{proof}
That $\Fscr_{3D}(\sigma)$ admits a unique minimizer was established in \cite{MP-KRM} under appropriate conditions on the constraint $n$. Furthermore, with the notations 
$$
\sigma_x=\TR_{\perp}(\sigma), \qquad \sigma_\perp=\TR_{x}(\sigma),
$$
for the partial traces w.r.t. $L^2(\Omega)$ and $L^2(0,L)$, respectively, and for any density operator $\sigma$ on $L^2((0,L) \times \Omega)$, the subaddivity of $- S(\varrho)$ yields, see \cite{Araki1970},
$$
S(\sigma)\geq S_x(\sigma_x)+S_\perp(\sigma_\perp).
$$
Hence, for any density operator $\sigma$  on $L^2((0,L) \times \Omega)$ with $n_\sigma=n$,
\bee
    \Fscr_{3D}(\sigma)&\geq& S_x(\sigma_x)+S_\perp(\sigma_\perp)+ \TR_x(\mathscr{H}_x \sigma_x)+\TR_\perp(\mathscr{H}_\perp \sigma_\perp)= \Fscr_x(\sigma_x)+\Fscr_\perp(\sigma_\perp)\\
&\geq & \Fscr_x(\varrho_\star)+\Fscr_\perp(\varrho_\perp).
\eee
Above, $\varrho_\star$ denotes $\varrho_\star[n_0]$ for simplicity, and we used that $n_\sigma=n$ implies $n_{\sigma_x}=n_0$. Finally, since  a direct calculation shows that $\Fscr_{3D}(\varrho_\star \otimes \varrho_\perp)=\Fscr_x(\varrho_\star)+\Fscr_\perp(\varrho_\perp)$, it follows that $\Fscr_{3D}(\sigma) \geq \Fscr_{3D}(\varrho_\star \otimes \varrho_\perp)$ for any density operator $\sigma$ satisfying the constraint $n_\sigma=n$. Since the eigenfunctions of $\mathscr{H}_\perp$ are complex exponentials as a consequence of the periodic boundary conditions, it follows that $n_{\varrho_\perp}=|\Omega|^{-1}$, and therefore that $n_{\varrho_\star \otimes \varrho_\perp}=n$. Hence, $\varrho_\star \otimes \varrho_\perp$ is the unique minimizer of $\Fscr_{3D}$ under the local constraint $n$.
\end{proof}

\bigskip

We are now in position to conclude. We need the following assumptions on the solutions to \fref{liou3D}: we suppose that (i) \fref{liou3D} admits a unique solution under appropriate conditions on the initial condition $\varrho_0$, and (ii) that this solution is obtained as the limit in proper sense as $k \to \infty$ of the sequence $\{\varrho_k\}_{k \in \Nm}$, that satisfies the linear problem
\begin{equation} \label{lin}
    i  \partial_t\varrho_{k+1} = [\mathscr{H}_{3D},\varrho_{k+1}] + i  (\varrho_e[\varrho_k]-\varrho_{k+1}), \qquad \varrho_{k+1}(t=0)=\varrho_0,
\end{equation}
Items (i) and (ii) are established in 1D in \cite{MP-JDE} without the uniqueness result, the latter being proven in Section \ref{convanal} in the present paper. The 3D case is still open.

We proceed by induction to obtain that $\varrho_{k+1}=\varrho_{x,k+1} \otimes \varrho_\perp$ where $\varrho_{x,k+1}$ verifies
\begin{equation*} 
    i  \partial_t\varrho_{x,k+1} = [\mathscr{H}_{x},\varrho_{x,k+1}] + i  (\varrho_\star[n_{\varrho_{x,k}}]-\varrho_{x,k+1}), \qquad \varrho_{x,k+1}(t=0)=\varrho_{x,0}.
\end{equation*}
For $k=0$, we have $n_{\varrho_0}=n_{\varrho_{x,0}} |\Omega|^{-1}$ with $\varrho_0=\varrho_{0,x} \otimes \varrho_\perp$, and therefore, according to Lemma \ref{subadd}, $\varrho_{e}[\varrho_0]=\varrho_\star[n_{\varrho_{x,0}}] \otimes \varrho_\perp$. Since \fref{lin} is linear and admits a unique solution, it follows that $\varrho_{1}=\varrho_{x,1} \otimes \varrho_\perp$ for an appropriate $\varrho_{x,1}$. Since the same reasoning applies for any $k$, we obtain that $\varrho_{k+1}=\varrho_{x,k+1} \otimes \varrho_\perp$. Using assumption (ii), it follows that the 3D solution $\varrho$ reads $\varrho(t)=\varrho_x(t) \otimes \varrho_\perp$, where $\varrho_x$ verifies the 1D equation
\begin{equation*} 
    i  \partial_t\varrho_{x} = [\mathscr{H}_{x},\varrho_{x}] + i  (\varrho_\star[n_{\varrho_{x}}]-\varrho_{x}), \qquad \varrho_{x}(t=0)=\varrho_{x,0}.
\end{equation*}

Note that we considered a linear potential $W(x)$ in this section, but the same approach holds for the 3D Poisson potential $V_{3D}$ since the resolution of the Laplace equation $\alpha^2 \Delta V_{3D}=n_\varrho=|\Omega|^{-1} n_{\varrho_x}$ with Dirichlet boundary conditions on $(0,L)$ and periodic on $\Omega$ yields $V_{3D}(x,y,z)=V_{3D}(x)$. 

This ends the justification of the 1D model.

\subsection{Proof of Lemma \ref{splitHbnd}} \label{proofsplitHbnd}
Given $\varrho_s^k$, we recall that one iteration of the splitting scheme reads
$$
\varrho_{s}(t) = U(t/2)W(t)U(t/2)\varrho^k_s, \qquad t \in [0,h],
$$
where $U(t)\sigma = e^{-i\Ham_0 t}\sigma e^{i\Ham_0 t}$, and $\varrho_1(t):=W(t) \sigma$ is the solution to
\be \partial_t \varrho_1 = \rhoeq{\varrho_1}-\varrho_1, \qquad \varrho_1(0)=\sigma. \label{colproof} \ee
\paragraph{Existence and uniqueness.} We show first that the $\varrho_s(t)$ above is well-defined and unique. We proceed iteratively. First, if $\varrho_s^0$ is a density operator in $\calE$, then so is $U(t/2) \varrho_s^0$ for all $t \geq 0$ since $U(t)$ preserves self-adjointness and positivity, and 
$$
\| U(t) \varrho_s^0 \|_{\calE}=\|\varrho_s^0\|_{\calE}.
$$
Let $\sigma:=U(h/2) \varrho_s^0$. Considering the collision subproblem $\varrho_1(t)=W(t) \sigma$, we recall that \fref{colproof} preserves the local density, and a consequence the equation is linear and admits as solution
\be \label{solcol}
\varrho_{1}(t) = (1-e^{-t})\rhoeq{\sigma} + e^{-t}\sigma,\ee
provided $\rhoeq{\sigma}$ exists and is unique. According to \cite[Theorem 2.1]{MP-JSP}, the latter holds when $\sigma \in \calE$, and when $n[\sigma](x)>0$ for all $x \in [0,1]$, yielding a unique $\rhoeq{\sigma} \in \calE_+$. We already know that $\sigma \in \calE_+$ from the previous step, and need to prove the lower bound. Following the assumptions of Theorem \ref{exist}, for any $t \geq 0$,
$$
n[U(t) \varrho_s^0]=n[f(\Ham_0)]+n[U(t) \delta \varrho] \geq \underline{n}+n[U(t) \delta \varrho] \qquad \textrm{(since $e^{ it \Ham_0}$ commutes with $f(\Ham_0)$)}, $$
and, under again the assumptions of Theorem \ref{exist}, we have
$$
\|n[U(t) \delta \varrho]\|_{L^\infty} \leq \|n[U(t)\delta \varrho]\|_{W^{1,1}} \leq 2 \|U(t) \delta \varrho\|_\calE = 2 \|\delta \varrho\|_\calE \leq \underline{n}/2.
$$
This shows that $n[U(t)\varrho_s^0](x) \geq \underline{n}/2$ for all $t$, and therefore that $\rhoeq{\sigma}$ exists in $\calE_+$ and is unique. Hence, $\varrho_1$ is well-defined, and as a consequence so is $\varrho_s^1=\varrho_s(h)$ in $\calE_+$. We now iterate over $k$. Since $\varrho_e[U(t)\varrho_s^0]$ is nonnegative, we have from \fref{solcol}, for all $\tau \geq 0$,
$$
n[U(\tau) \varrho_{s}^1 ] \geq e^{-h} n[U(\tau+h) \varrho_s^0] \geq e^{-h} \underline{n}/2,
$$
which allows us to construct $\varrho_e[ U(h/2) \varrho_s^1] \in \calE_+$ and therefore $\varrho_s^2$. Iterating, we then find $\varrho_s^k \in \calE$ and, from the version of \fref{solcol} at step $k$, 
\be \label{lb2}
n[U(\tau)\varrho_{s}^k ] \geq e^{-h} n[U(\tau+h) \varrho_s^{k-1}] \geq e^{- k h} n[U(\tau+k h) \varrho_s^{0}] \geq   e^{- k h} \underline{n}/2 \geq e^{- T} \underline{n}/2,
\ee
which proves the lower bound on $n[U(\tau)\varrho_s^k]$ for all $k$ and all $\tau \geq 0$. We have therefore obtained a unique solution to the splitting scheme in $\calE$ satisfying the lower bound announced in the lemma.
\paragraph{Uniform bounds.} We derive now a bound in $\calH$ that is uniform in $k$ and $h$. For this, we need first uniform bounds in $\calJ_1$ and in $\calE$. The one in $\calJ_1$ is direct as $U(t)$ is an isometry in $\calJ_1$ and \fref{colproof} preserves trace, and therefore
$$
\| \varrho_s^k \|_{\calJ_1}=\|\varrho_s^0\|_{\calJ_1}.
$$
For the bound in $\calE$, we remark that $U(t)$ is an isometry in $\calE$, and that we have the following bound from Proposition 2.2 in \cite{MP-JDE}:
$$
\| \varrho_e[\sigma] \|_\calE \leq C \|\sqrt{n[\sigma]}\|^2_{H^1} \leq C+C \|\sigma\|_{\calE},
$$
where $C$ is independent of $\sigma$. With the above definition of $\varrho_1$, this yields,
%\[\begin{split} \| \varrho_{1}(t)\|_\calE &\leq (1-e^{-t})\|\rhoeq{\sigma}\|_\calE + e^{-t}\|\sigma\|_\calE \leq C' t(1 +\|\sigma\|_\calE),
 % \end{split}\]
\[\begin{split} \| \varrho_{1}(t)\|_\calE &\leq (1-e^{-t})\|\rhoeq{\sigma}\|_\calE + e^{-t}\|\sigma\|_\calE \leq C t (1+\|\sigma\|_\calE) + e^{-t}\| \sigma \|_\calE,
  \end{split}\]
since $1-e^{-t} \leq t$ for $t\geq 0$. Going back to the splitting solution $\varrho_s$, we therefore obtain
$$
\| \varrho_s^{k+1} \|_\calE \leq Ch+C h \|\varrho_s^{k}\|_\calE + e^{-h}\| \varrho_s^{k} \|_\calE \leq Ch+e^{Ch} \| \varrho_s^{k} \|_\calE,
$$
%$$
%\| \varrho_s^{k+1} \|_\calE \leq C'h (1+\|\varrho_s^{k}\|_\calE).
%$$
Iterating, it follows that
%$$
%\| \varrho_s^{N_T} \|_\calE \leq C T e^{C T } \| \varrho^0 \|_\calE,
%$$
$$
\| \varrho_s^{N_T} \|_\calE \leq Ch \sum_{k=0}^{N_T-1} e^{Ck h} \| \varrho^0 \|_\calE \leq C T e^{C T } \| \varrho^0 \|_\calE,
$$
which provides us with a uniform bound in $\calE$. We move on now to the $\calH$ bound, and use the following result from \cite{MP-JDE}: let $\sigma \in \calH$, with $\|\sigma\|_\calE \leq \alpha_0$ and $n[\sigma] \geq \alpha_1>0$. Then,
\be \label{sublin}
\| \varrho_e[\sigma] \|_\calH \leq C_{\alpha_0,\alpha_1}\|\sigma\|_{\calH}.
\ee
With the above definition of $\varrho_1$, this yields
\[\begin{split} \| \varrho_{1}(t)\|_\calH &\leq (1-e^{-t})\|\rhoeq{\sigma}\|_\calH + e^{-t}\|\sigma\|_\calH \leq C t \|\sigma\|_\calH + e^{-t}\| \sigma \|_\calH,
  \end{split}\]
where the constant $C$ is independent of $k$ and $h$ since the lower bound in \fref{lb2} and the bound in $\calE$ are uniform in $k$ and $h$. Going back to the splitting solution $\varrho_s$, we therefore obtain
$$
\| \varrho_s^{k+1} \|_\calH \leq C h \|\varrho_s^{k}\|_\calH + e^{-h}\| \varrho_s^{k} \|_\calH \leq e^{Ch} \| \varrho_s^{k} \|_\calH,
$$
 Iterating, it follows that
$$
\| \varrho_s^{k} \|_\calH \leq e^{C k h } \| \varrho^0 \|_\calH.
$$
This ends the proof.
\subsection{Proof of Lemma \ref{localerror}} \label{prooflocalerror}
Before proceeding with the proof, the following generalized Gronwall Lemma will be useful. The proof of the general result can be found in \cite{gronwall}. 

\begin{lemma}[Gronwall]
Let $f:[0,T]\rightarrow \mathbb{R}$ be continuous and satisfy the  inequality, 
\[f(t) \leq M + \int_0^t e^{-(t-s)}(f(s))^{\gamma}ds, \qquad \gamma \in (0,1),\]
where $M\geq 0$. Then, the following estimate holds
\[ f(t) \leq \Phi^{-1}\Big(\Phi(M) + 1-e^{-t}\Big),\]
where $\Phi(u) = \frac{1}{1-\gamma}u^{1-\gamma}$ and $\Phi^{-1}(w) = (1-\gamma)^\frac{1}{1-\gamma}w^\frac{1}{1-\gamma}$. 
\end{lemma}

The following two Lemmas can be found in \cite{MP-JDE} and will be used in the proof.

\begin{lemma}[Lemma 6.4 in \cite{MP-JDE}] \label{lem64}
Let $\varrho \in \calH$, self-adjoint and nonnegative. Then, 
\[ \jnorm{U(t)\varrho -\varrho} \leq Ct\hnorm{\varrho} \quad \text{for all } t\geq 0.\]
\end{lemma}

\medskip

The result below shows that the map $\varrho \mapsto \varrho_e[\varrho]$ is at least of H\"older regularity $1/8$ in $\calJ_2$.
\begin{lemma}[Corollary 5.8 in \cite{MP-JDE}] \label{cor58}

Let $\varrho_{1}$ and $\varrho_{2}$ be two density operators in $\calH$. Let $M_{0}\in(0,\infty)$ be such that 
\[\hnorm{\varrho_{1}}+\hnorm{\varrho_{2}} \leq M_{0}, \quad \text{ and   } M_{0}^{-1}\leq n_{\varrho_i}, \quad \text{ for all } x\in[0,1], \quad i=1,2.\]
Then,
$$
\|\varrho_e[\varrho_1]-\varrho_e[\varrho_2]\|_{\calJ_2} \leq C \|\varrho_1-\varrho_2\|^{1/8}_{\calJ_2},
$$
where $C$ is independent of $\varrho_1$ and $\varrho_2$.
\end{lemma}

\medskip

We can now proceed with the proof. According to \fref{intsplit} and \fref{intex}, the error $e_k(t):=\varrho(t_k+t)-\varrho_s(t_k+t)$ for $t \in [0,h]$, with the notation $e_k:=e_k(0)$, verifies
\bee
e_{k}(t)&=& e^{-t} U(t) (\varrho^{k}-\varrho^{k}_{s})\\
&&+ \int_{0}^{t} e^{-(t-u)} \left(U(t-u)\rhoeq{\varrho(t_k+u)}-U(t/2)\rhoeq{U(t/2)\varrho^k_s} \right)du,
\eee
where $\varrho^k=\varrho(t_k)$ and $\varrho_s^k=\varrho_s(t_k)$. Taking the $\calJ_2$ norm and using the fact that $U(t)$ is an isometry on $\calJ_2$, we find for $t \in [0,h]$,
%\begin{equation*}
 % \begin{split}
    \bee
    \hsnorm{e_{k}(t)} &\leq& e^{-t}\hsnorm{\varrho^{k}-\varrho^{k}_{s}}\\
    &&+\int_{0}^{t} e^{-(t-u)}\hsnorm{U(t-u)\rhoeq{\varrho(t_k+u)}-U(t/2)\rhoeq{U(t/2)\varrho^k_s}}du\\
    & \leq &e^{-t}\hsnorm{e_{k}}\\
    &&+ \int_{0}^{t} e^{-(t-u)}\hsnorm{U(t-u)\rhoeq{\varrho(t_k+u)}-U(t-u)\rhoeq{U(t/2)\varrho^{k}_{s}}}du \\
    &&+\int_{0}^{t}e^{-(t-u)}\hsnorm{U(t-u)\rhoeq{U(t/2)\varrho^{k}_{s}}-U(t/2)\rhoeq{U(t/2)\varrho^k_{s}}}du\\[2mm]
      &=:& e^{-t} \hsnorm{e_{k}} + I_{1}(t) + I_{2}(t).
\eee
  %  \end{split}
 %\end{equation*}
First, consider the integral given by $I_{2}$. We have, since $U$ is an isometry on $\calJ_2$,
\begin{equation*}
\begin{split}
    I_{2}(t) &= \int_{0}^{t}e^{-(t-u)}\hsnorm{U(t/2-u)\rhoeq{U(t/2)\varrho^{k}_{s}}-\rhoeq{U(t/2)\varrho^{k}_{s}}}du\\
    &\leq C\int_{0}^{t}e^{-(t-u)}|t/2-u|\hnorm{\rhoeq{U(t/2)\varrho^{k}_{s}}}du\\
    &\leq C\int_{0}^{t}e^{-(t-u)}|t/2-u|\hnorm{\varrho^{k}_{s}}du\\
& \leq Ct^{2}\hnorm{\varrho^0}.
\end{split}
\end{equation*}
The first inequality is thanks to Lemma \ref{lem64} and the fact that $\calJ_2 \subset \calJ_1$. The second inequality is due to the sublinear estimate $\hnorm{\rhoeq{U(t/2)\varrho^{k}_{s}}} \leq C\hnorm{U(t/2)\varrho^{k}_{s}}$ stated in \fref{sublin}, which holds provided $n[U(t/2)\varrho^{k}_{s}] \geq \alpha>0$ and $U(t/2)\varrho^{k}_{s}$ is bounded uniformly in $\calE$. These two facts are obtained in Lemma \ref{splitHbnd} as $\calH \subset \calE$. The last inequality is due to estimate \fref{eq:unfHbnd} in Lemma \ref{splitHbnd}.

Now, consider the integral term $I_{1}$. We apply Lemma \ref{cor58} as both $\varrho$ and $U(t/s)\varrho^k_s$ belong to $\calH$ and their respective local densities are uniformly bounded from below according to Theorem \ref{exist} and Lemma \ref{splitHbnd}. Then, with $\gamma=1/8$,
\begin{equation*}
\begin{split}
    I_{1}(t) &= \int_{0}^{t}e^{-(t-u)} \hsnorm{\rhoeq{\varrho(t_k+u)}-\rhoeq{U(t/2)\varrho^k_{s}}}du\\
    &\leq C \int_{0}^{t}e^{-(t-u)}\hsnorm{\varrho(t_k+u)-U(t/2)\varrho^k_s}^{\gamma}du\\
    &\leq C \int_{0}^{t}e^{-(t-u)}\hsnorm{\varrho(t_k+u)-\varrho_{s}(t_k+u)}^{\gamma}du \\
&\qquad + C\int_{0}^{t} e^{-(t-u)}\hsnorm{\varrho_{s}(t_k+u)-U(t/2)\varrho^k_{s}}^{\gamma}du\\
    &=: T_{1}(t) + T_{2}(t).
\end{split}
\end{equation*}
The term $T_{1}$ will be handled further with the Gronwall Lemma. For $T_{2}$, we remark first that from \fref{intsplit} and Lemma \ref{splitHbnd},
$$
\|\varrho_{s}(t_k+u)-U(u)\varrho_s^k\|_{\calJ_2} \leq C h, \qquad \forall u \in [0,h]. 
$$
Then, using again Lemma \ref{splitHbnd} and Lemma \ref{lem64}, we find, for $t \in [0,h]$,
\bee
        T_{2}(t) &\leq& C h^{1+\gamma}+\int_{0}^{t} e^{-(t-u)}\hsnorm{\varrho_s^k-U(t/2-u)\varrho^k_{s}}^{\gamma}du \\
        &\leq& C h^{1+\gamma}+C\int_{0}^{t} e^{-(t-u)}  |t/2-u|^\gamma \hnorm{\varrho^k_{s}}^\gamma du \\
&\leq& Ch^{1+\gamma}.
 \eee
Collecting all estimates, we have for $t \in [0,h]$,
\begin{equation*}
    \begin{split}
        \hsnorm{e_{k}(t)}&\leq e^{-h}\hsnorm{e_{k}} + C h^{1+\gamma}  + Ch^{2}+\int_{0}^{t}e^{-(t-u)}\hsnorm{e_{k}(s)}^{\gamma}du\\
        &=: M_{k,h} + \int_{0}^{t}e^{-(t-u)}\hsnorm{e_{k}(s)}^{\gamma}du.
    \end{split}
  \end{equation*}
The generalized Gronwall Lemma then yields, using that $(x+y)^\beta \leq C_\beta (x^\beta+y^\beta)$ for $x,y \geq 0$ and $\beta \geq 1$, for $t \in [0,h]$,
  \begin{equation*}
    \begin{split}
        \hsnorm{e_{k}(t)}
                &\leq (1-\gamma)^{\frac{1}{1-\gamma}}\Big(\frac{1}{1-\gamma}M_{k,h}^{1-\gamma}+1-e^{-t}\Big)^{\frac{1}{1-\gamma}} \\
        & \leq M_{k,h} + Ch^{\frac{1}{1-\gamma}}= e^{-h}\hsnorm{e_{k}}+ C (h^{1+\gamma}+h^2+h^{\frac{1}{1-\gamma}}).
    \end{split}
\end{equation*}
This ends the proof. 

\subsection{Semi-classical approximation} \label{semiC}

We obtain relation \fref{semiN} by using pseudo-differential calculus, and need for this to extend the problem to the whole $\Rm$. We remain at a formal level. Let then $\chi_\beta$ be a smooth function over $\Rm$ such that $0 \leq \chi_\beta \leq 1$, with $\chi_\beta(x)=0$ for $x \notin [0,1]$, and $\chi_\beta(x)=1$ for $x \in [\beta^\gamma, 1-\beta^\gamma]$, for some $\gamma<1$. Denoting by $o(1)$ quantities that are negligible in appropriate sense when $\beta \ll 1$, we have, for any smooth function $\varphi$,
\be \label{eqtra}
\int_0^1 n[e^{-(\Ham_0+\Ascr)}](x) \varphi(x)dx= \Tr(e^{-(\Ham_0+\Ascr)} \varphi)=\Tr( \chi_\beta e^{-(\Ham_0+\Ascr))} \chi_\beta \varphi)+o(1),
\ee
since the size of the support of $1-\chi$ is less than $\beta^\gamma$. Next, we remark that the function $\exp(-(\Ham_0+\Ascr)) \chi_\beta \varphi$ is equal to the function $u(t=1,x)$, with $u$ solution to
$$
\partial_t u=-(\Ham_0+\Ascr)u, \qquad u(t=0,x)=\chi_\beta(x) \varphi(x).
$$
Consider then the operator $\beta^2 \Delta -\Ascr_e$ defined on $\Rm$, with $\Ascr_e=\Ascr$ on $[0,1]$ and $\Ascr_e(x)=x^2$ for $x \notin [0,1]$. With $v=\chi_\beta u$, we find that
$$
\partial_t v=-(\Ham_0+\Ascr)v+R_\beta=(\beta^2 \Delta -\Ascr_e)v+R_\beta, \qquad v(t=0,x)=\chi^2_\beta(x) \varphi(x),
$$
where
$$
R_\beta=-\beta^2\,u \Delta \chi_\beta -2\beta^2 \nabla \chi_\beta \nabla u,
$$
and the hypotheses on $\chi_\beta$ yield $R_\beta=o(1)$. This shows that
$$
\chi_\beta e^{-(\Ham_0+\Ascr)} \chi_\beta \varphi=e^{(\beta^2 \Delta -\Ascr_e)} \chi^2_\beta \varphi+o(1),
$$
and as a consequence, with \fref{eqtra},
\be \label{appN}
\int_0^1 n[e^{-(\Ham_0+\Ascr)}](x) \varphi(x)dx=\int_0^1 n[e^{(\beta^2 \Delta -\Ascr_e)}](x) \chi^2_\beta(x) \varphi(x)dx+o(1).
\ee
We are now in position to use pseudo-differental calculus and find an approximation for $n[e^{(\beta^2 \Delta -\Ascr_e)}]$, which is well-defined since $e^{(\beta^2 \Delta -\Ascr_e)}$ is trace class because of the confining potential $\Ascr_e$. It is shown in \cite{QET} that
$$
n[e^{(\beta^2 \Delta -\Ascr_e)}](x)=\frac{1}{\sqrt{4 \pi} \beta} e^{-\Ascr_e}(x)+o(1), \qquad \forall x \in \Rm,
$$
and therefore, according to \fref{appN},
$$
n[e^{-(\Ham_0+\Ascr)}](x)=n[e^{(\beta^2 \Delta -\Ascr_e)}](x) \chi^2_\beta(x)+o(1)=\frac{1}{\sqrt{4 \pi} \beta} e^{-\Ascr}(x)\chi^2_\beta(x)+o(1).
$$
This gives \fref{semiN}. Regarding the approximate functional, we set $\varphi=1$ and find
\bee
\int_0^1 n[e^{-(\Ham_0+\Ascr)}](x) dx&=&\frac{1}{\sqrt{4 \pi} \beta} \int_0^1 e^{-\Ascr}(x)\chi^2_\beta(x) dx+o(1)\\
&=&\frac{1}{\sqrt{4 \pi} \beta} \int_0^1 e^{-\Ascr}(x)dx+o(1).
\eee
The functional $J_{\rm{approx}}$ is finally obtained by spatial discretization. This ends this section.

\end{appendix}
\bibliographystyle{plain}
\bibliography{../../bibliography.bib}
%\printbibliography
\end{document}